\documentclass[12pt,a4paper]{article}
\usepackage{xcolor}
\usepackage{times}
\usepackage{calc}
\usepackage{graphicx} 
\usepackage{subfigure}
\usepackage{amssymb}
\usepackage{amsmath}
\usepackage{float}
\usepackage{latexsym}
\usepackage{psfrag}
\usepackage{dsfont}
\usepackage{rotating}
\usepackage{longtable}
\usepackage[numbers]{natbib}

\usepackage [english]{babel}
\usepackage{hyperref}

\usepackage[a4paper,left=3cm,right=2cm,top=2.5cm,bottom=2.5cm]{geometry}

\newcommand{\DG}{{\mathrm{DG}}}
\newcommand{\FV}{{\mathrm{FV}}}

\newcommand{\mat}[1]{\underline{{#1}}}

\def\d{\mathrm{d}}

\usepackage{accents}
\newcommand\acclrvec[1]{\accentset{\,\leftrightarrow}{#1}}	

\newcommand{\state}[1]{\mathbf{#1}}
\newcommand{\blockvec}[1]{\acclrvec{{\mathbf #1}}\ }		

\newcommand{\bigpartialderiv}[2]{ \frac{\partial {#1}}{\partial {#2} } }
\newcommand\norm[1]{\left\lVert#1\right\rVert}

\newcommand{\numfluxb}[1]{\hat{\mathbf{#1}} }

\title{A Subcell Finite Volume Positivity-Preserving Limiter for DGSEM Discretizations of the Euler Equations}

\author{Andrés M. Rueda-Ramírez
\thanks{Department of Mathematics and Computer Science, University of Cologne, Weyertal 86-90, 50931 Cologne, Germany.}
\thanks{Corresponding author: Andrés M. Rueda-Ramírez (aruedara@uni-koeln.de).}
 \and Gregor J. Gassner $^*$
	 }

\date{}

\begin{document}

\maketitle

\begin{abstract}
In this paper, we present a positivity-preserving limiter for nodal Discontinuous Galerkin disctretizations of the compressible Euler equations.
We use a Legendre-Gauss-Lobatto (LGL) Discontinuous Galerkin Spectral Element Method (DGSEM) and blend it locally with a consistent LGL-subcell Finite Volume (FV) discretization using a hybrid FV/DGSEM scheme that was recently proposed for entropy stable shock capturing.
We show that our strategy is able to ensure robust simulations with positive density and pressure when using the standard and the split-form DGSEM.
Furthermore, we show the applicability of our FV positivity limiter in extremely under-resolved vortex dominated simulations and in problems with shocks.
\end{abstract}

\section{Introduction}

Discontinuous Galerkin (DG) methods 
provide high-order accuracy and local conservation on arbitrary grids when dealing with advection-dominated problems.
The combination of high-order DG methods with high-order explicit Runge-Kutta schemes for the time integration leads to a highly accurate and easily parallelizable family of methods, e.g.,  \cite{Cockburn2001,Hindenlang2012}.
However, it has been observed that in extremely under-resolved fluid flow simulations, or in the presence of shocks, the oscillations of the high-order polynomials can lead to non-physical negative values of density or pressure and cause the scheme to crash. 

A particular class of DG methods, known as Discontinuous Galerkin Spectral Element Methods (DGSEM), can be combined with a split-form representation of the governing equations \cite{Gassner2016} or entropy-stable fluxes \cite{Carpenter2014} to obtain non-linear stability and hence increased robustness.
Nevertheless, these schemes still suffer from positivity issues in the presence of shocks \cite{Hennemann2020,Rueda-Ramirez2020} and under-resolved turbulence.

An effective and rather commonly used strategy to preserve physical (positive) values of density and pressure in high-order DG methods was proposed by \citet{Zhang2010}.
This method shrinks the high-order numerical solution in each element around its mean value until a valid solution is obtained.
Although this method is often very effective, we observed that it sometimes causes large jumps in the numerical solution, which may trigger instabilities and lead to restrictions in the time-step size.

In this paper, we propose a novel positivity-preserving limiter for DGSEM discretizations that is based on a subcell FV discretization.
We make use of the hybrid FV/DGSEM scheme that was recently developed by \citet{Hennemann2020} for shock capturing.
The hybrid FV/DGSEM scheme blends a Legendre-Gauss-Lobatto (LGL) DGSEM discretization with a compatible LGL-subcell Finite Volume (FV) scheme in an element-wise manner. 
In the present work, we present a simple strategy to compute the blending coefficient of each element at the end of every Runge-Kutta stage, such that the positivity condition is never violated.


\section{The Euler Equations of Gas Dynamics} \label{sec:Euler}

The Euler equations can be written as a system of conservation laws
\begin{equation} \label{eq:ConsLaw}
\bigpartialderiv{\state{u}}{t} + \nabla \cdot \blockvec{f}(\state{u}) = \state{0},
\end{equation}
where the state quantities are density, momentum and total energy, $\state{u} = [\rho, \rho \vec{v}, \rho E]^T$, and the flux reads
\begin{equation}
\blockvec{f} (\mathbf{u})=
\begin{bmatrix}
\rho \vec{v} \\
\rho \vec{v} \otimes \vec{v} + \mat{I} p \\
\vec{v} (\rho E + p)
\end{bmatrix},
\end{equation}
where $\mat{I}$ is the identity matrix, and the pressure is computed with the calorically perfect gas assumption,
%
$
p=(\gamma-1) \rho e,
$
$\gamma$ is the heat capacity ratio, and $e = E - \norm{\vec{v}}^2/2$ is the internal energy.

\section{Numerical Discretization Scheme} \label{sec:Discretization}

\subsection{Spatial Discretization: Hybrid FV/DGSEM Discretization}

We make use of the blended scheme by  \citet{Hennemann2020} that combines a high-order DGSEM discretization with a first-order FV method on the corresponding LGL-subcell grid to obtain an accurate scheme with robust shock-capturing properties. As both schemes are constructed on the same subcell LGL distribution, they are directly compatible and thus the semi-discrete version of \eqref{eq:ConsLaw} is given by
\begin{equation} \label{eq:ODEblended}
\dot{\state{u}} = 
(1-\alpha) \dot{\state{u}}^{\DG} 
+ \alpha \dot{\state{u}}^{\FV}, 
\end{equation}
where $\dot{\state{u}}^{\DG}$ is the high-order DG discretization of the PDE, $\dot{\state{u}}^{\FV}$ is the first-order FV discretization, and $\alpha \in [0,1]$ is an element-local blending coefficient that depends on a shock indicator.
In the present work, we are not necessarily interested in using this blending approach for shock capturing, but, as we will show, it can be used for positivity preservation.


To obtain the DGSEM-discretization of \eqref{eq:ConsLaw}, the simulation domain is tessellated into $K$ non-overlapping elements and all variables are approximated within each element by piece-wise Lagrange interpolating polynomials of degree $N$ on LGL nodes.
These polynomials are continuous in each element and may be discontinuous at the element interfaces.
Furthermore, \eqref{eq:ConsLaw}  is multiplied by an arbitrary polynomial (test function) of degree $N$ and numerically integrated by parts inside each element of the mesh with an LGL quadrature rule of $N+1$ points on a reference element, $\xi \in [-1,1]$, to obtain
\begin{equation} \label{eq:ODEDG}
\dot{\state{u}}_j^{\DG} = 
\state{F}^{\mathrm{Vol}}_j
+ \frac{\delta_{jN}}{J \omega_j} \left( \state{f}_N - \numfluxb{f}_{(N,R)} \right)
- \frac{\delta_{j0}}{J \omega_j} \left( \state{f}_0 - \numfluxb{f}_{(L,0)} \right),\quad j=0,...,N.
\end{equation}
In \eqref{eq:ODEDG}, $\omega_j$ is the reference-space quadrature weight, $\delta_{ij}$ is the Kronecker delta, $J$ is the geometry mapping Jacobian from reference space to physical space, $\numfluxb{f}_{(i,j)} = \numfluxb{f}(\state{u}_i,\state{u}_j)$ is the surface numerical flux, which we evaluate with the left and right outer states, $L$ and $R$, to account for the jump of the solution across cell interfaces, and $\state{F}_{\mathrm{Vol}}$ is the volume integral term.

From \eqref{eq:ODEDG}, we can recover the standard and the split-form (flux differencing) DGSEM discretizations by choosing the volume integral term as
\begin{equation}
\state{F}^{\mathrm{Vol,Std}}_j = 
-\frac{1}{J} \sum_{i=0}^N D_{ji} \state{f}_i, 
~~~~ \mathrm{or} ~~~
\state{F}^{\mathrm{Vol,Split}}_j = 
-\frac{2}{J} \sum_{i=0}^N D_{ji} \state{f}^{*}_{(j,i)},
\end{equation}
respectively, where $D_{ji}= \ell'_i(\xi_j)$ is the so-called derivative matrix, defined in terms of the Lagrange interpolating polynomials, $\{ \ell_i \}_{i=0}^N$, and $\state{f}^{*}_{(j,i)} = \state{f}^{*}(\state{u}_j,\state{u}_i)$ is a volume numerical two-point numerical flux.
We can construct a \textit{provably} entropy-stable (ES) DGSEM scheme by choosing an entropy conserving flux for $\state{f}^{*}$, e.g. \cite{Chandrashekar2013}, and an entropy stable flux for $\numfluxb{f}$, e.g. \cite{Gassner2016}.

The compatible LGL co-located first-order FV discretization proposed by \citet{Hennemann2020}, which interprets the nodal values as the subcell mean values, reads
\begin{equation} \label{eq:ODEFV}
\dot{\state{u}}_j^{\FV} =
\frac{1}{J \omega_j} \left(
  \numfluxb{f}_{(j,j-1)}
- \numfluxb{f}_{(j,j+1)}
\right),
\end{equation}
where the numerical fluxes on the element boundaries match with the DGSEM boundary values and are evaluated with the left and right outer states,
%
$\numfluxb{f}_{(0,-1)} := \numfluxb{f}(\state{u}_0,\state{u}_L)$, 
$\numfluxb{f}_{(N,N+1)} := \numfluxb{f}(\state{u}_N,\state{u}_R)$.

\subsection{Time-Integration Scheme}

We obtain the evolution of the state variables over time by applying a high-order explicit strong stability-preserving Runge-Kutta (SSP-RK) method to \eqref{eq:ODEblended}.

In the SSP-RK methods of \citet{Spiteri2002}, the solution is updated in every RK stage, $s$, with 
\begin{align} \label{eq:NextU}
\state{u}^{s+1} &=
a_{ss} \state{u}^s + \Delta t^s \dot{\state{u}}^s +
\sum_{i=1}^{s-1} \left( a_{si} \state{u}^i + \Delta t b_{si} \dot{\state{u}}^i \right),
\end{align}
where $a_{si}$ and $b_{si}$ are the RK coefficients, $\Delta t^s = b_{ss} \Delta t$ is the so-called time-step size of the $s$ RK stage.

\section{Positivity-Preserving Scheme} \label{sec:PositivityMethod}

First-order FV schemes are positivity preserving in density and pressure for the Euler equations if a forward Euler time-discretization is used and an appropriate Riemann solver is selected \cite{Perthame1996}.
Therefore, if $\alpha=1$, the solution that is obtained is safe (positive in density and pressure),
\begin{equation} \label{eq:NextUsafe}
\state{u}_{\rm{safe}}^{s+1} = 
a_{ss} \state{u}^s + \Delta t^s \dot{\state{u}}^{s,\FV} +
\sum_{i=1}^{s-1} \left( a_{si} \state{u}^i + \Delta t b_{si} \dot{\state{u}}^i \right).
\end{equation}

As a result, it is possible to modify $\alpha$ locally, such that the discrete solution remains positive in density and pressure throughout the simulation.
We propose a goal condition based on the computed safe solution
\begin{equation} \label{eq:rhop_restriction}
\rho \ge \beta \rho_{\rm{safe}}, 
~~~~
p \ge \beta p_{\rm{safe}},
\end{equation}
with $\beta \in (0, 1]$. 
In other words, we allow only a certain deviation (controlled by the value $\beta$) below the safe solution, which is a stricter requirement than positivity. 
We do not consider upper bounds for density or pressure, as we are focusing on positivity. 
Note that upper bounds are equally straightforward to impose and that upper bounds might be necessary to avoid oscillatory solutions.

\subsection{Density Correction}

Let us obtain the corrected blending coefficient for density first. 
This is easy, because the density is one of the conserved quantities, and as such, it depends linearly on the blending coefficient.
Using \eqref{eq:NextU} and \eqref{eq:NextUsafe}, we get for every degree of freedom of an element,
\begin{equation} \label{eq:rho_safeDif}
\rho_{\rm{safe}} - \rho = \Delta t^s (\dot{\rho}^{s,\FV} - \dot{\rho}^{s}).
\end{equation}
Summing and subtracting $\beta \rho_{\rm{safe}}$ we obtain
%
\begin{equation} \label{eq:rho_safeDif2}
\rho_{\rm{safe}} - \rho_{\rm{new}} +  a_{\rho} = \Delta t^s (\dot{\rho}^{s,\FV} - \dot{\rho}^{s}),
\end{equation}
where $\rho_{\rm{new}}=\beta \rho_{\rm{safe}}$ and $a_{\rho} = \beta \rho_{\rm{safe}} - \rho$.

If $a_{\rho} > 0$, requirement \eqref{eq:rhop_restriction} is not fulfilled and $\alpha$ must be corrected such that $\rho = \rho_{\rm{new}}$.
We rewrite \eqref{eq:rho_safeDif2} as
\begin{equation}
\rho_{\rm{safe}} - \rho_{\rm{new}} = \Delta t^s (\dot{\rho}^{s,\FV} - \dot{\rho}^{s,\mathrm{new}}),
\end{equation}
where
%
$
\dot{\rho}^{s,\mathrm{new}} = \dot{\rho}^s + \frac{a_{\rho}}{\Delta t^s},
$
such that we can compute a new blending coefficient \textit{a posteriori} to ensure the fulfillment of the density criterion,
\begin{equation}
\alpha_{\rm{new}} = \alpha + \frac{a_{\rho}}{\Delta t^s (\dot{\rho}^{s,\FV} - \dot{\rho}^{s,\DG})},
\end{equation}
where $\alpha$ is the blending coefficient that was used to obtain $\dot{\state{u}}$.

The corrected solution and its time derivative can be then computed with $\Delta \alpha = \alpha_{\mathrm{new}} - \alpha$ as
\begin{equation} \label{eq:CorrectU_Ut}
\state{u}^{s+1}_{\mathrm{new}} = \state{u}^{s+1} + \Delta \alpha \Delta t^s (\dot{\state{u}}^{\FV} - \dot{\state{u}}^{\DG}),
~~~~
\dot{\state{u}}^{s+1}_{\mathrm{new}} = \dot{\state{u}}^{s+1} + \Delta \alpha (\dot{\state{u}}^{\FV} - \dot{\state{u}}^{\DG}).
\end{equation}

Note that, since we use element-local blending coefficients, $\alpha_{\mathrm{new}}$ must be taken as the maximum of the corrected blending coefficients computed for all the degrees of freedom of an element.

\subsection{Pressure Correction}
After correction of density, we check the resulting approximation for pressure. If condition \eqref{eq:rhop_restriction} is not fulfilled for the pressure, one must solve the following nonlinear equation for $\alpha_{\rm{new}}$,
\begin{equation}
g(\alpha_{\rm{new}}) = p(\state{u}_{\rm{new}}(\alpha_{\rm{new}})) - \beta p_{\rm{safe}} = 0,
\end{equation}
which can be done using an iterative Newton's method,
\begin{equation} \label{eq:NewtonAlpha}
\alpha^{n+1}_{\mathrm{new}} = \alpha^n -\frac{g(\alpha^n)}{\partial p (\alpha^n) / \partial \alpha}.
\end{equation}
where $\partial p / \partial \alpha$ can be easily obtained using the chain rule,
\begin{equation}
\bigpartialderiv {p}{\alpha} = \bigpartialderiv {p}{\state{u}} \bigpartialderiv {\state{u}} {\alpha}.
\end{equation}

The first term in the right-hand side is obtained from \eqref{eq:ODEblended} and \eqref{eq:NextU}, and the second term can be derived from the equation of state,
\begin{equation}
\bigpartialderiv {\state{u}} {\alpha} = \Delta t^s (\dot{\state{u}}^{\FV} - \dot{\state{u}}^{\DG}),
~~~~
\bigpartialderiv {p}{\state{u}} = (\gamma - 1) \left( \frac{1}{2} \norm{\vec{v}}^2,~ -\vec{v},~ 1 \right).
\end{equation}

We finally compute the corrected solution and its time derivative using  \eqref{eq:CorrectU_Ut} and the new element-wise blending coefficient, $\alpha_{\mathrm{new}}$, which is again taken as the maximum of the corrected blending coefficients for all the degrees of freedom of an element.

In our proposed \textit{a posteriori} positivity-preserving limiter, we first check density values and correct invalid states, and then check pressure values and correct if needed.
The reason is that, in many cases, performing the density correction removes invalid values of pressure automatically.
Furthermore, the density correction is computationally cheaper than the pressure correction since the former can be computed with a linear, single-step process, while the latter needs an iterative method.

\section{Numerical Results} \label{sec:Results}

For the examples of this section, we use the 2-register SSP-RK of fourth order and five stages of \citet{Spiteri2002}, where the time step is computed with the typical CFL condition for the DGSEM \cite{Krais2019} with CFL $=0.5$. 
To control the positivity conditions \eqref{eq:rhop_restriction}, we choose the value $\beta=0.1$.

\subsection{Kelvin-Helmholtz Instability}

To test the capabilities of the FV positivity limiter, we first apply it to simulate an inviscid Kelvin-Helmholtz instability (KHI).
This setup is quite challenging for a high-order DG method, as it is severely under-resolved with an effective Reynolds number Re $=\infty$.
The initial condition is given by
\begin{align} \label{eq:KHI_IniCond}
\rho (t=0) &= \frac{1}{2} 
+ \frac{3}{4} B,
&
p (t=0) &= 1,
\nonumber \\
v_1 (t=0) &= \frac{1}{2} \left( B-1 \right),
&
v_2 (t=0) &= \frac{1}{10} \sin(2 \pi x),
\end{align}
with
$
B=
\tanh \left( 15 y + 7.5 \right) - \tanh(15y-7.5).
$

The simulation domain, $\Omega=[-1,1]^2$, is complemented with periodic boundary conditions.
For this example, we tessellate $\Omega$ using $K=64^2$ quadrilateral elements, represent the solution with polynomials of degree $N=3$, and run the simulation until $t=25$.

The initial condition, \eqref{eq:KHI_IniCond}, provides a Mach number Ma $\le 0.6$, which makes compressibility effects relevant, but does not cause shocks to develop in the simulation domain. Thus, we do not apply shock capturing by connecting $\alpha$ to a troubled cell indicator, e.g. \cite{Hennemann2020}, but instead just use $\alpha = 0$ as our baseline high-order discretization.

The Euler equations are discretized using the standard DGSEM and the entropy-stable DGSEM (ES-DGSEM). In both frameworks, we employ a traditional Rusanov scheme for the surface numerical fluxes, $\numfluxb{f}$, and for the ES-DGSEM we use the entropy-conserving and kinetic energy preserving flux of \citet{Chandrashekar2013}  for $\state{f}^*$.

\begin{figure}[htb]
\centering
\subfigure[Entire simulation.]{
	\includegraphics[trim=20 0 60 0,clip,width=0.42\linewidth]{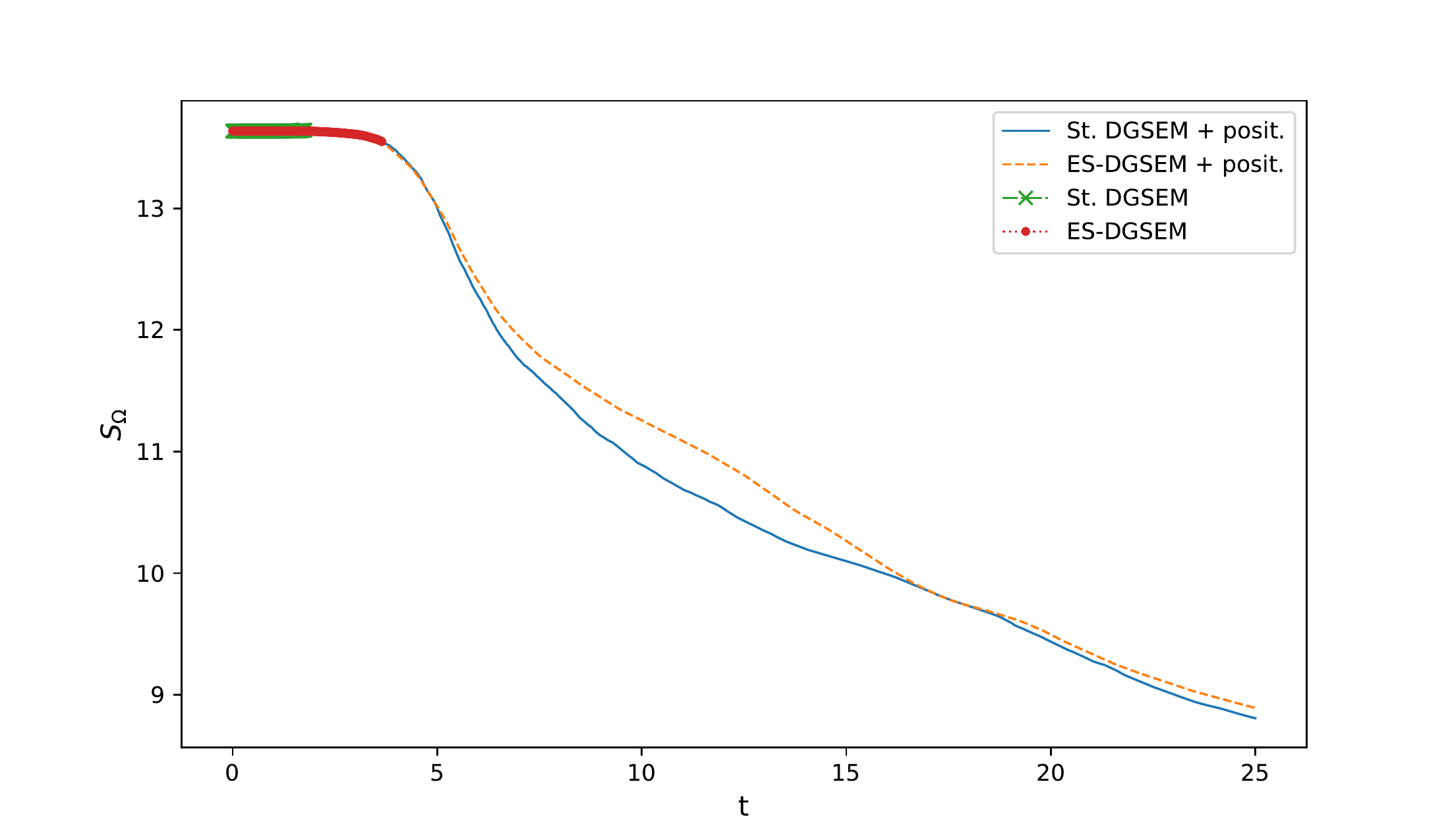}
}
\subfigure[Beginning of the simulation.]{
	\includegraphics[trim=20 0 60 0,clip,width=0.42\linewidth]{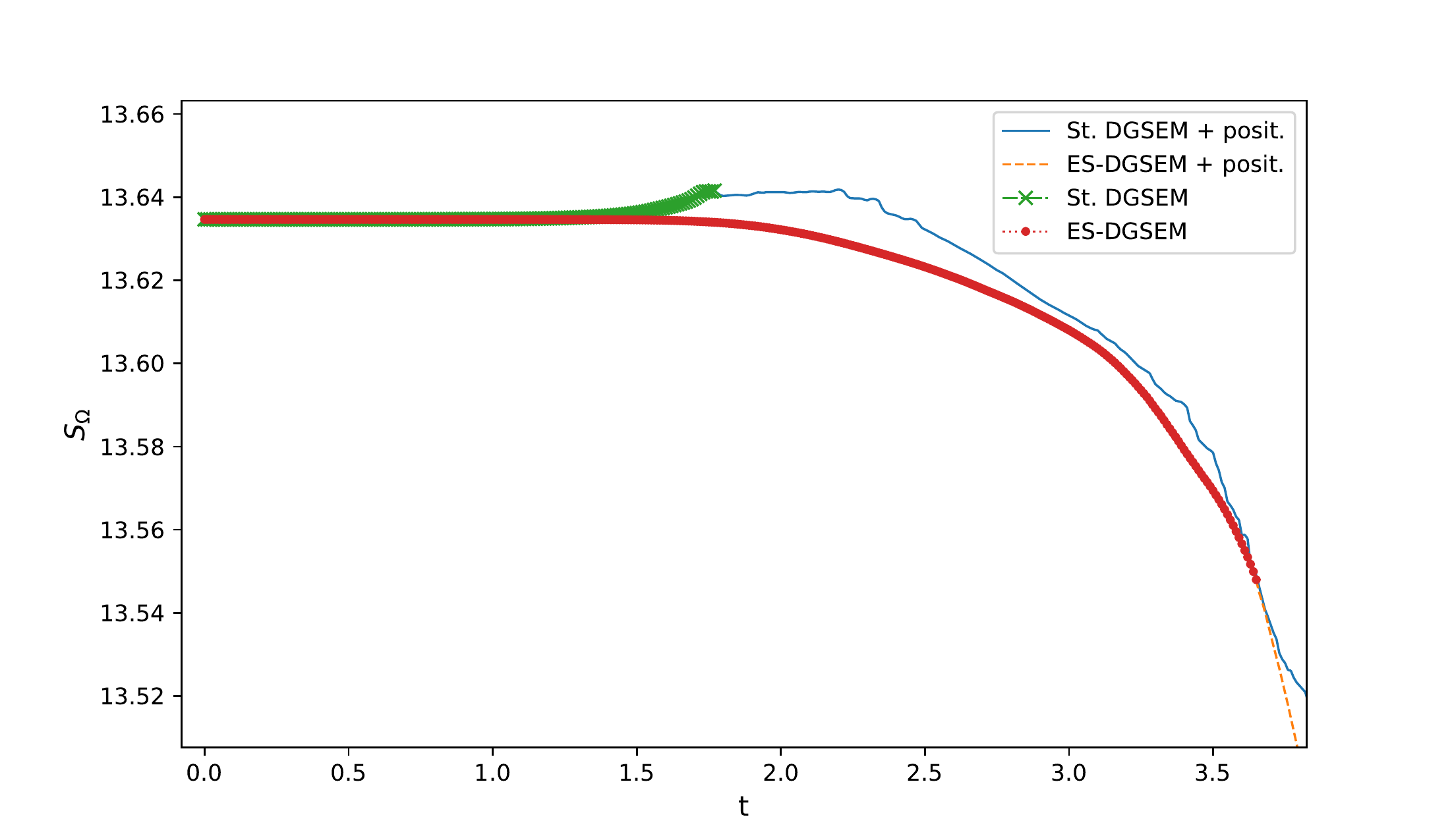}
}
\caption{Total entropy evolution over time (a) and detailed zoom at the beginning of the simulation (b) for the Kelvin-Helmholtz instability simulation.}
\label{fig:Entropy}
\end{figure}

Figure \ref{fig:Entropy} shows the evolution of the total mathematical entropy over time in the entire simulation domain,
\begin{equation} \label{eq:totalEntropy}
S_{\Omega}(t) = -\int_{\Omega} \frac{\rho s}{\gamma-1} \d x,
\end{equation}
where $s=\ln\left(p \rho^{-\gamma}\right)$ is the thermodynamic entropy.
Note that the conventional use of the minus sign in \eqref{eq:totalEntropy} implies that $S_{\Omega}$ decreases in the simulation domain if the second law of thermodynamics is fulfilled.

As can be observed in Figure \ref{fig:Entropy}, the ES-DGSEM fulfills the second law of thermodynamics, but the standard DGSEM does not. Note, however, that entropy stability is only guaranteed for positive solutions and that the ES-DGSEM has no inbuilt positivity preservation. It can be observed (close up plot (b) in the figure) that both schemes indeed crash at the beginning of the simulation due to positivity problems. However, with the proposed positivity limiter strategy, both schemes are able to run until the end of the simulation.

\begin{figure}[htb]
\centering
\subfigure[Maximum blending coefficient.]{
	\includegraphics[trim=30 0 60 0,clip,width=0.42\linewidth]{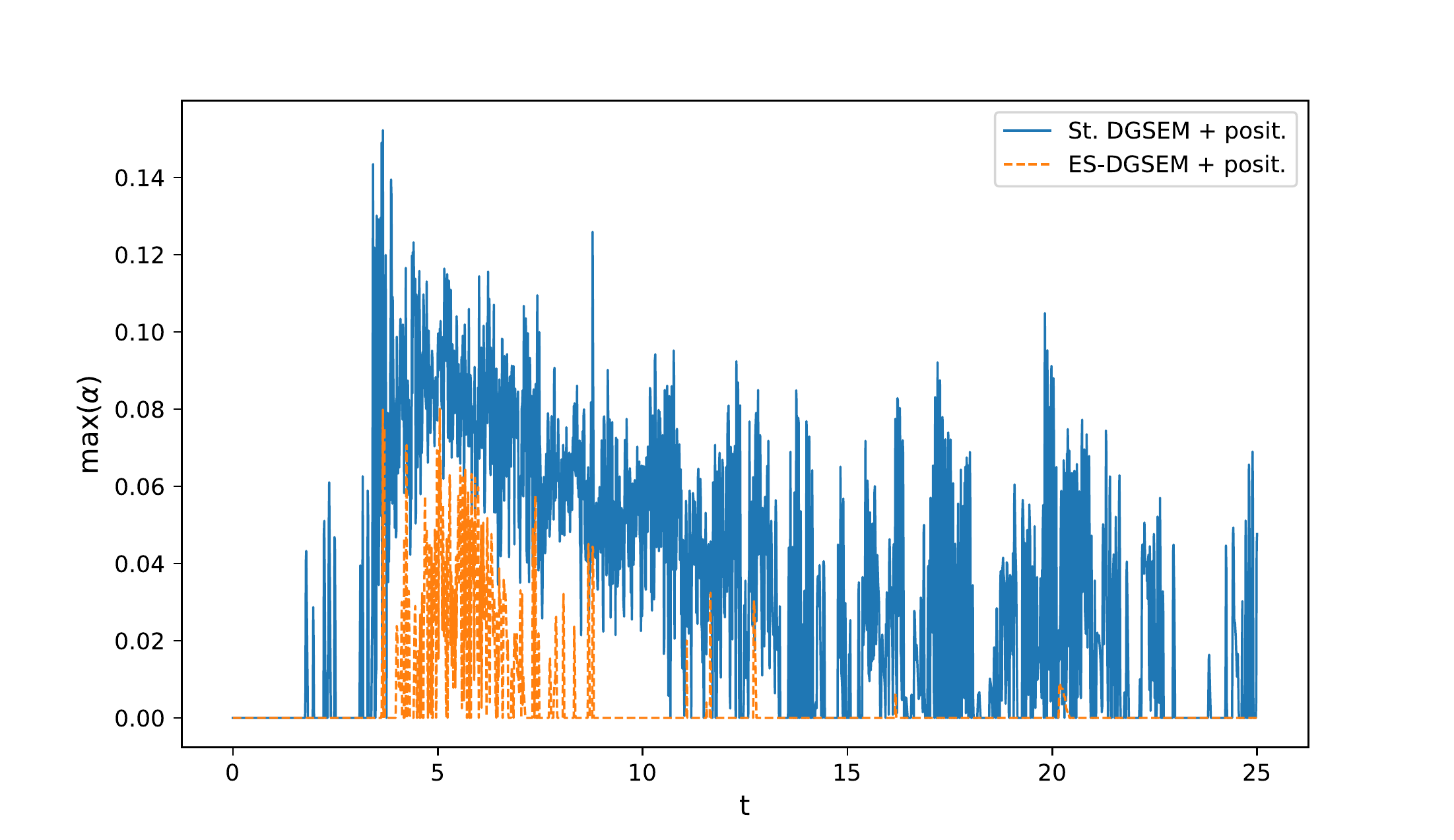}
}
\subfigure[Amount of low order method used in the domain.]{
	\includegraphics[trim=30 0 60 0,clip,width=0.42\linewidth]{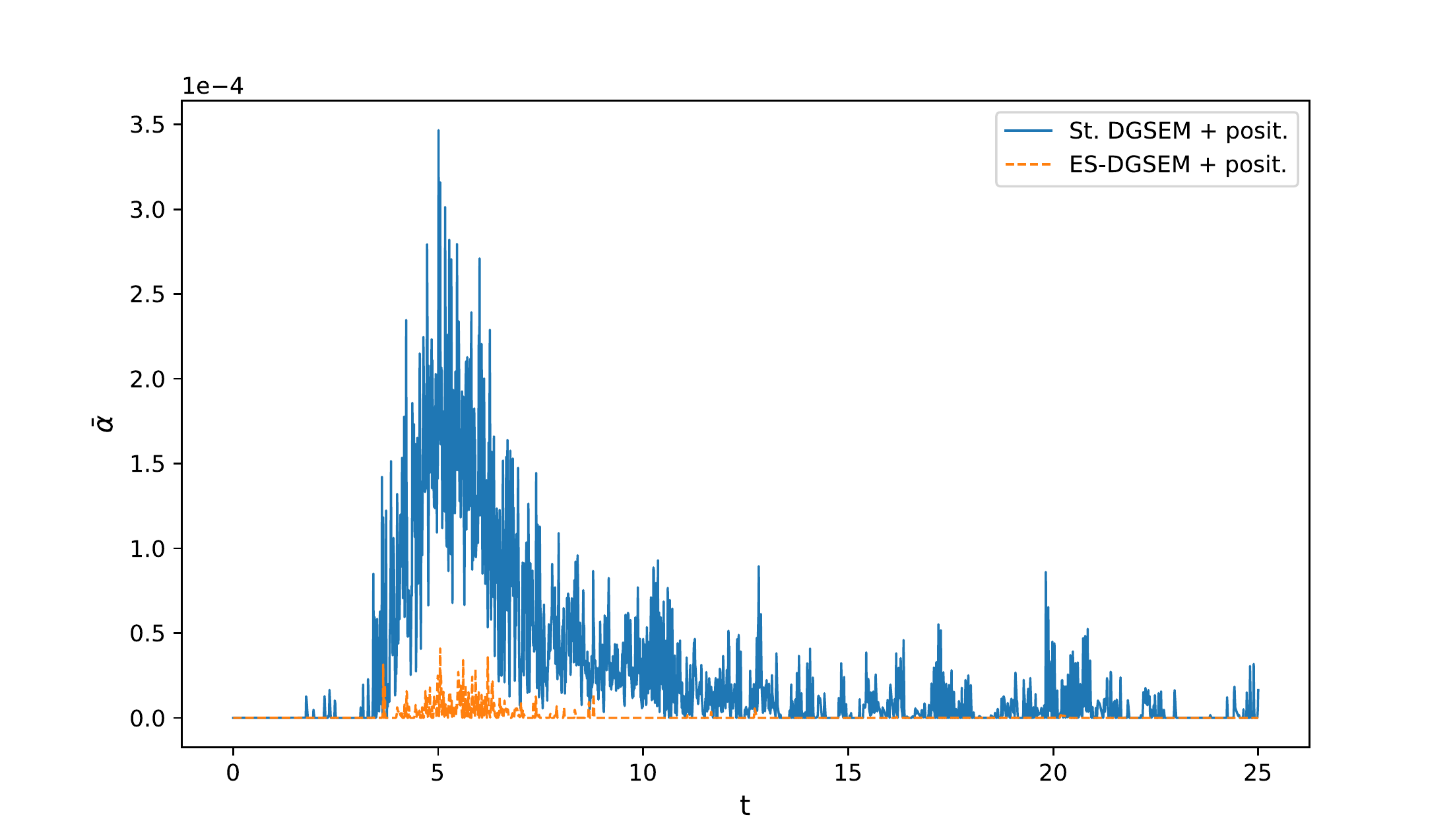}
}
\caption{Evolution of the blending coefficient over time for the Kelvin-Helmholtz instability simulation.}
\label{fig:BlendingCoeff}
\end{figure}

Figure \ref{fig:BlendingCoeff} illustrates the evolution of the blending coefficient over time in the whole domain using a sample time $\Delta \tau = 0.01$.
The quantities $\max (\alpha)$ and $\bar \alpha$ are defined as
\begin{equation} \label{eq:maxAlpha}
\max (\alpha) (t) = \max_{\tau \in [t-\Delta \tau, t]} \left( \max_{k=1}^K \alpha_k(\tau) \right),
~~~~
\bar \alpha (t) = \frac{1}{n_s} \sum_{s=1}^{n_s} \left( \frac{1}{K} \sum_{k=1}^K \alpha^s_k \right),
\end{equation}
where $k$ denotes the element index, $\alpha^s_k$ is the blending coefficient of element $k$ at the RK stage $s$, and $n_s$ is the number of RK stages taken from $t-\Delta \tau$ to $t$.
A value $\bar \alpha = 1$ means that the entire domain uses a first-order FV method, whereas $\bar \alpha = 0$ means that it uses only the high-order DG method.

In Figure \ref{fig:BlendingCoeff}(a) it is evident that only a small amount of first-order dissipation is needed to stabilize the simulation, as the maximum blending coefficient for any element of the domain is always below $\alpha \le 0.15$ for the standard DGSEM, and $\alpha \le 0.08$ for the ES-DGSEM.

%

The FV stabilization of the high-order schemes is added very locally.
As can be observed in Figure \ref{fig:BlendingCoeff}(b), the amount of the domain that is discretized with a first-order FV method is less than $0.035 \%$ during the entire simulation when the base-line scheme is the standard DGSEM, and less than $0.005 \%$ when the base-line scheme is the ES-DGSEM.
Moreover, in the latter case, the entire domain is discretized with a pure ES-DGSEM method during most of the simulation.

\begin{figure}[ht]
\centering
\includegraphics[trim=700 1250 700 90,clip,width=0.4\linewidth]{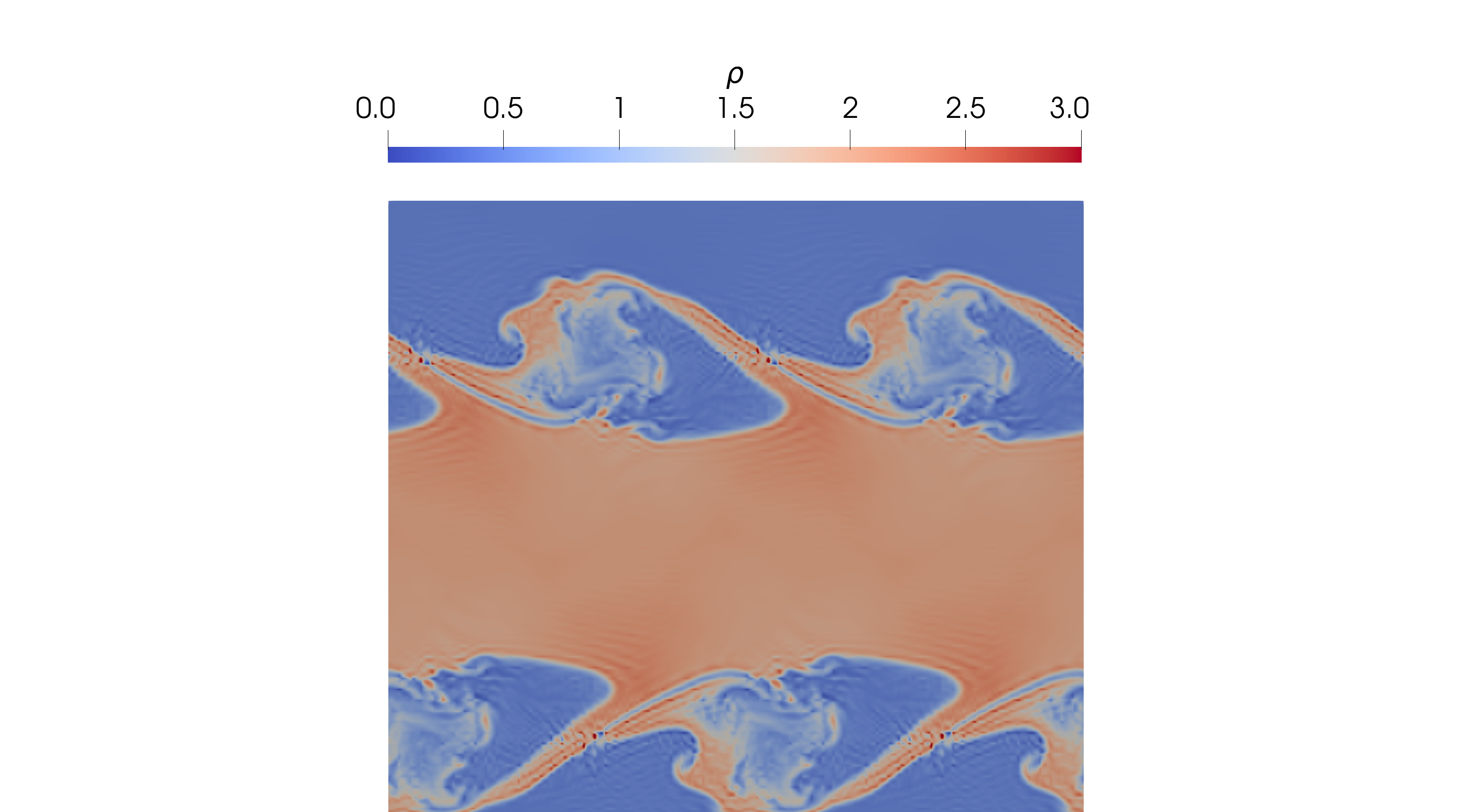}
\includegraphics[trim=700 1250 700 90,clip,width=0.4\linewidth]{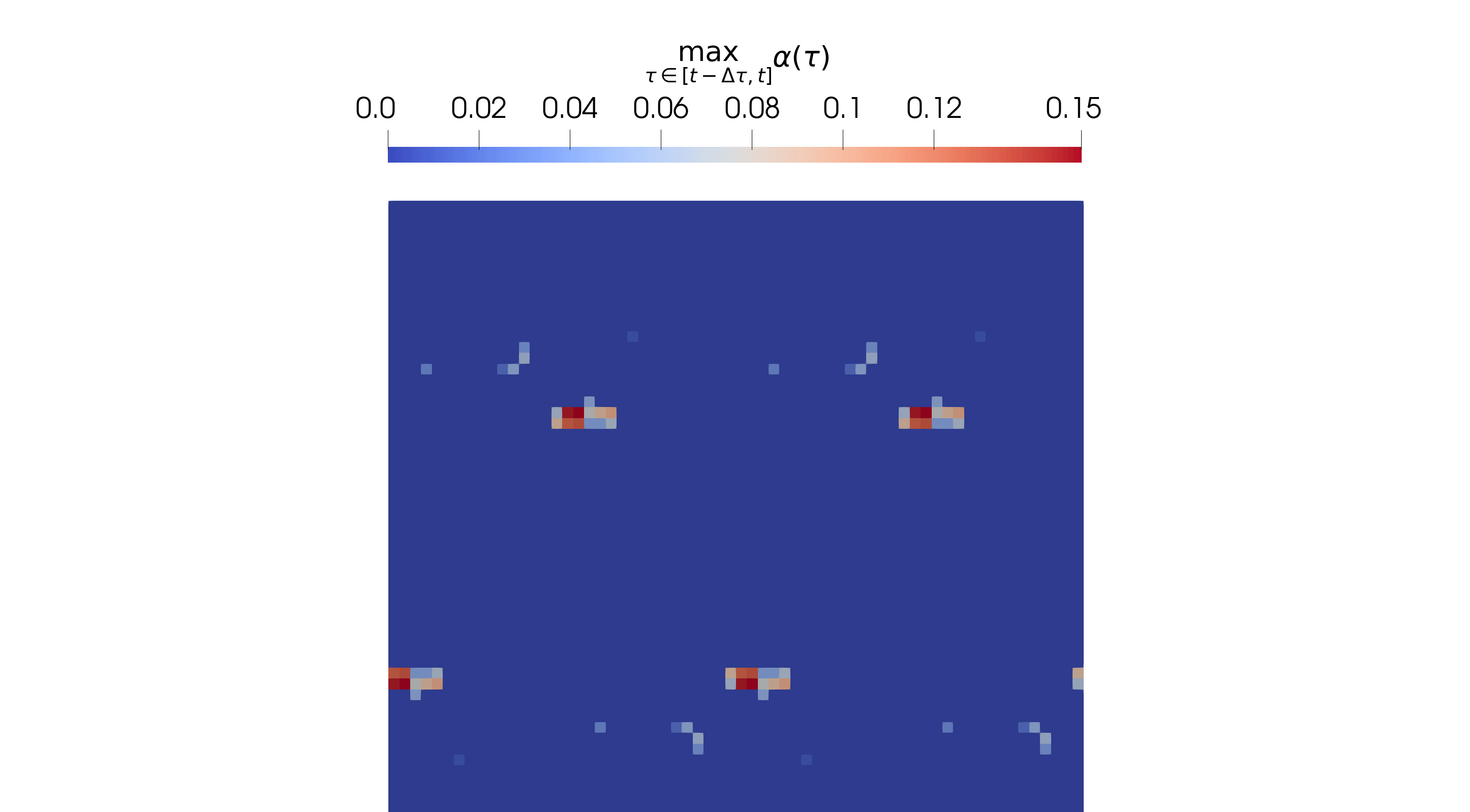}
\subfigure[St. DGSEM - Density.]{
	\includegraphics[trim=600 100 600 100,clip,width=0.4\linewidth]{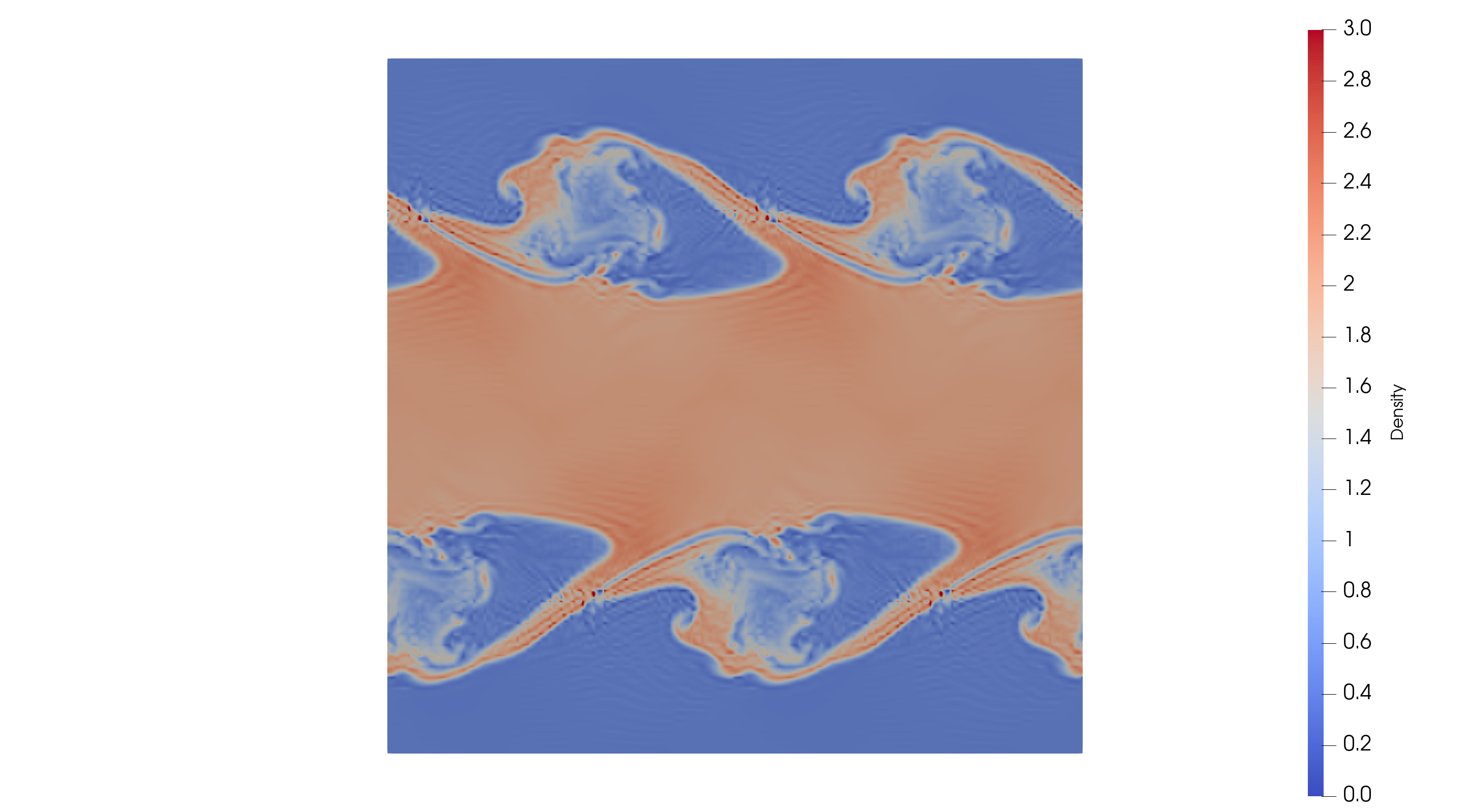}
}
\subfigure[St. DGSEM - Blending coefficient.]{
	\includegraphics[trim=600 100 600 100,clip,width=0.4\linewidth]{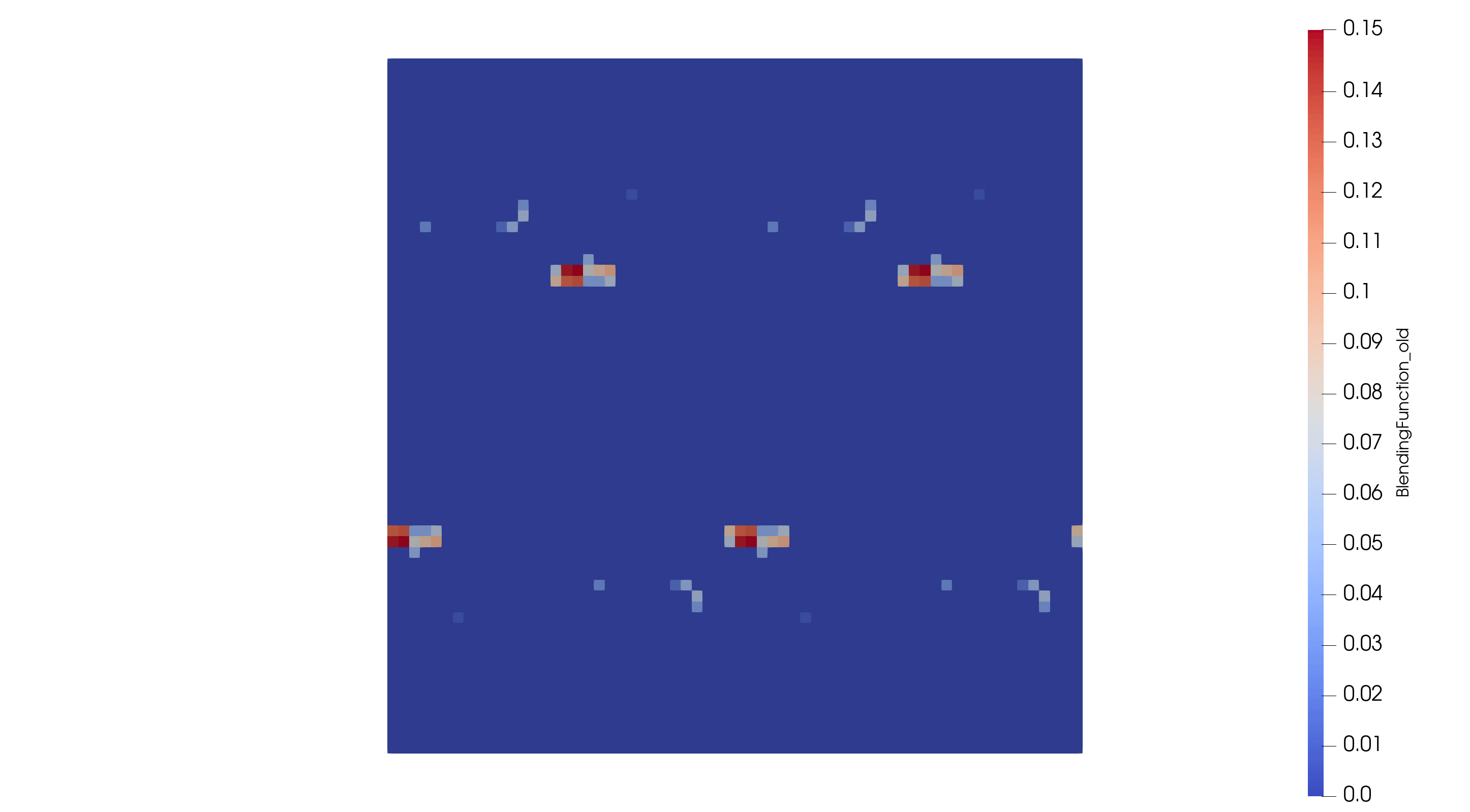}
}
\subfigure[ES-DGSEM - Density.]{
	\includegraphics[trim=600 100 600 100,clip,width=0.4\linewidth]{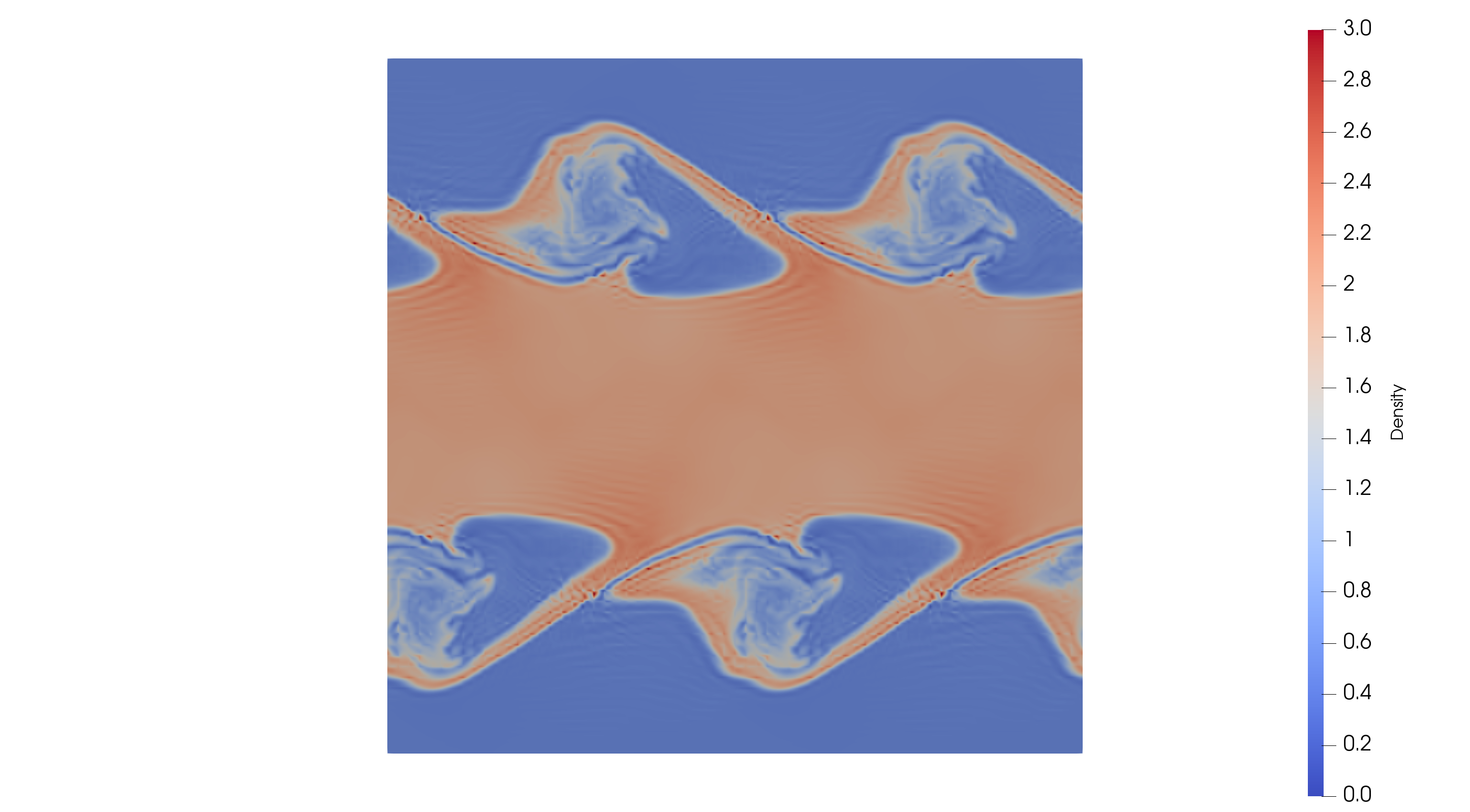}
}
\subfigure[ES-DGSEM - Blending coefficient.]{
	\includegraphics[trim=600 100 600 100,clip,width=0.4\linewidth]{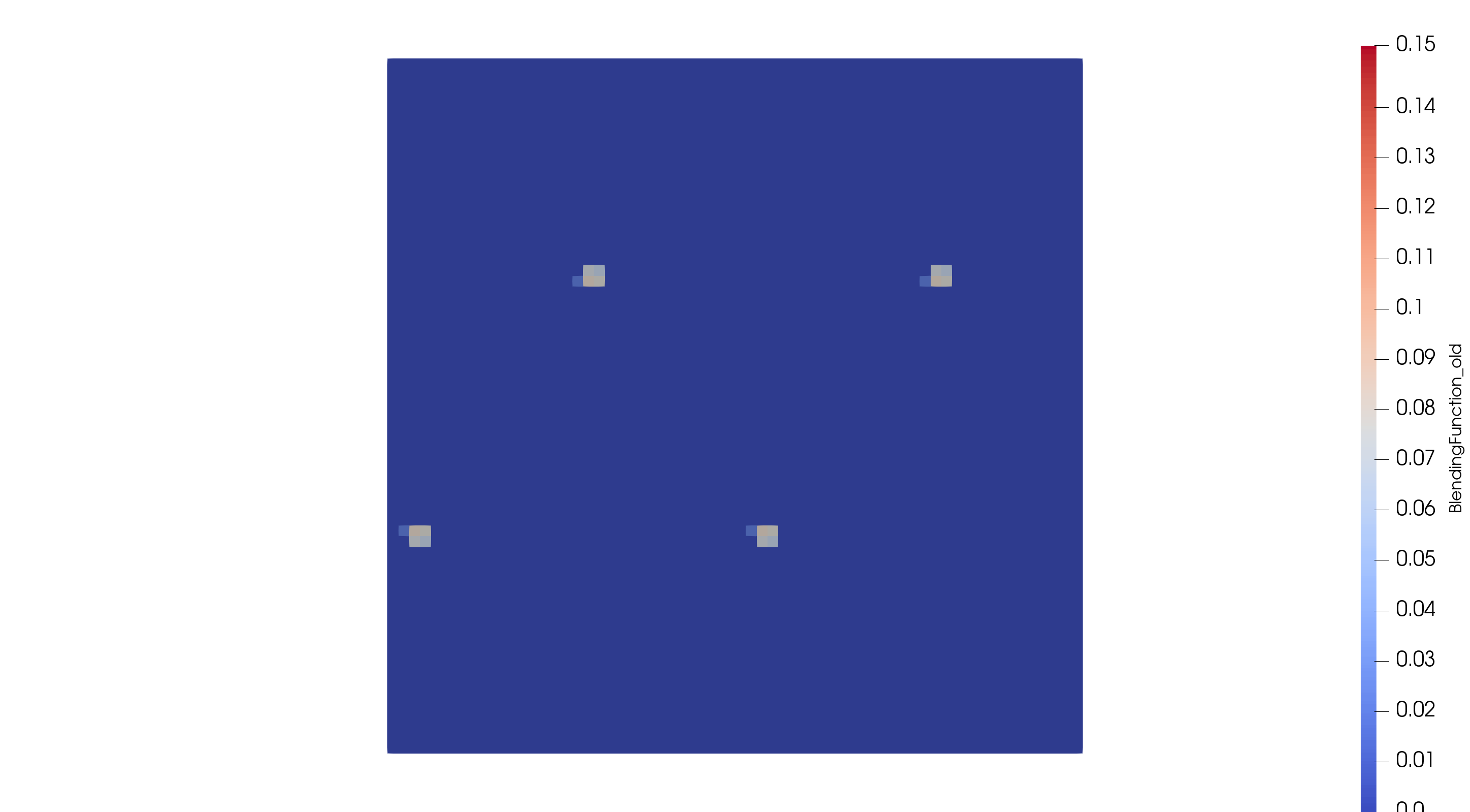}
}
\caption{Snapshot of the KHI simulations at $t=3.7$ for the standard DGSEM (top) and the ES-DGSEM (bottom). 
The blending coefficient is taken as the maximum over all the SSP-RK stages of a sampling interval $\Delta \tau = 0.01$.}
\label{fig:KHI_Contours}
\end{figure}

Finally, Figure \ref{fig:KHI_Contours} shows the density and the blending coefficient of both stabilized high-order schemes for a snapshot at $t=3.7$.
Since the blending occurs sporadically and in small amounts, it is very difficult to capture a blending coefficient $\alpha > 0$ at a given snapshot.
Therefore, we plot the maximum blending coefficient of each element of the domain over a sampling time $\Delta \tau=0.01$.
It can be observed in Figure \ref{fig:KHI_Contours} that the positivity limiter preserves the symmetries of the problem, and that the standard DGSEM needs more stabilization and is more distorted than the ES-DGSEM.

\subsection{Sedov Blast}

The Sedov blast problem describes the evolution of a radially symmetrical blast wave that expands from an initial concentration of density and pressure into a homogeneous medium.
For the initial condition, we assume a gas in rest, $v_1 (t=0) = v_2 (t=0) =0$, with a Gaussian distribution of density and pressure,
\begin{align} \label{eq:Sedov_IniCond}
\rho (t=0) &= \rho_0 + \frac{1}{4 \pi \sigma_{\rho}^2} \exp \left( -\frac{1}{2}  \frac{r^2}{\sigma_{\rho}^2} \right)
&
p (t=0) &= p_0 + \frac{\gamma - 1}{4 \pi \sigma_p^2} \exp \left( -\frac{1}{2}  \frac{r^2}{\sigma_p^2} \right),
\end{align}
where we choose $\sigma_{\rho}=0.25$ and $\sigma_p=0.15$.
Furthermore, the ambient density is set to $\rho_0=1$ and the ambient pressure to $p_0=10^{-5}$.

We complement the simulation domain, $\Omega=[-1.5,1.5]^2$, with periodic boundary conditions.
Furthermore, we tessellate $\Omega$ using $K=64^2$ quadrilateral elements, represent the solution with polynomials of degree $N=7$, and run the simulation until $t=20$.

The smooth initial condition defined by \eqref{eq:Sedov_IniCond} quickly evolves into an expanding shock front that bounces back into the domain because of the periodic boundary conditions. 
This generates a complex shock pattern with multiple shock-shock interactions, which causes the transition to turbulence (Re $=\infty$).

To deal with the strong shocks of this simulation, we use the indicator proposed by \citet{Hennemann2020}, which is a modification of \cite{Persson2006}.
We use the gas pressure for the indicator, such that $\alpha$ reads
\begin{equation} \label{eq:ShockIndicator}
\alpha = \frac{1}{1+\exp \left( \frac{-s}{\mathbb{T}}(\mathbb{E}-\mathbb{T})\right)},
\end{equation}
where $s=9.21024$ is the so-called sharpness, $\mathbb{T}(N)=0.5 \cdot 10^{-1.8 (N+1)^{0.25}}$ is the threshold, and $\mathbb{E}$ is the total energy of the two highest modes of the pressure solution.
%
Furthermore, we perform one spatial propagation sweep of the blending coefficient  to avoid large jumps in the blending coefficient from element to element.
The blending coefficient of element $k$ is adjusted with
%
$
\alpha_k  = \max_E \{ \alpha, 0.5 \alpha_E \},
$
where $\alpha_E$ denotes the blending coefficient of any neighbor element to $k$.

In this example, we use the HLLE scheme \cite{Einfeldt1988} for the surface numerical fluxes of both frameworks, $\numfluxb{f}$, which has been shown to be positivity preserving \cite{Einfeldt1991}.
For the ES-DGSEM we use the entropy-conserving and kinetic energy preserving flux of \citet{Chandrashekar2013}  for $\state{f}^*$.

\begin{figure}[ht]
\centering
\subfigure[Entire simulation.]{
	\includegraphics[trim=15 0 40 0,clip,width=0.42\linewidth]{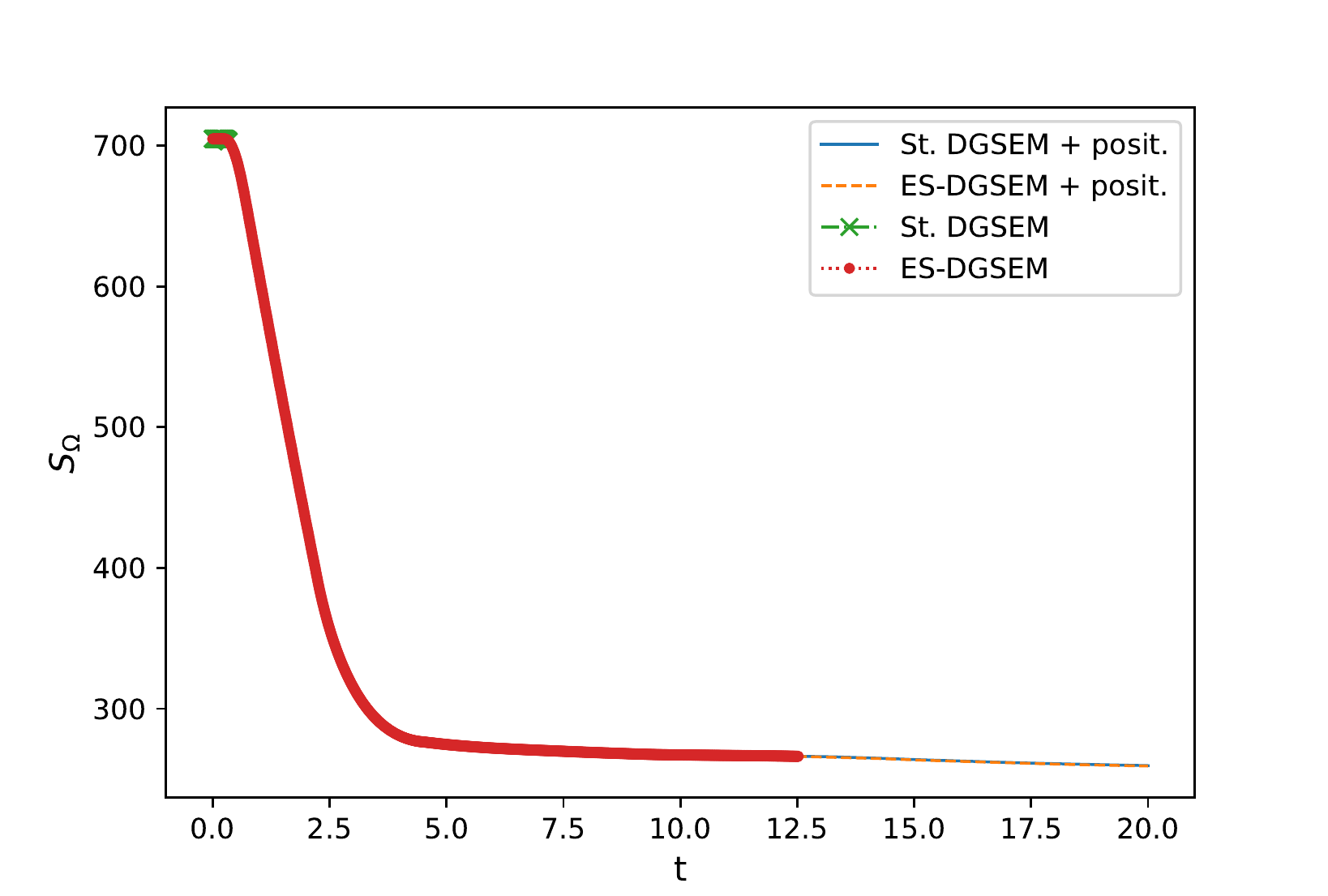}
}
\subfigure[Beginning of the simulation.]{
	\includegraphics[trim=15 0 40 0,clip,width=0.42\linewidth]{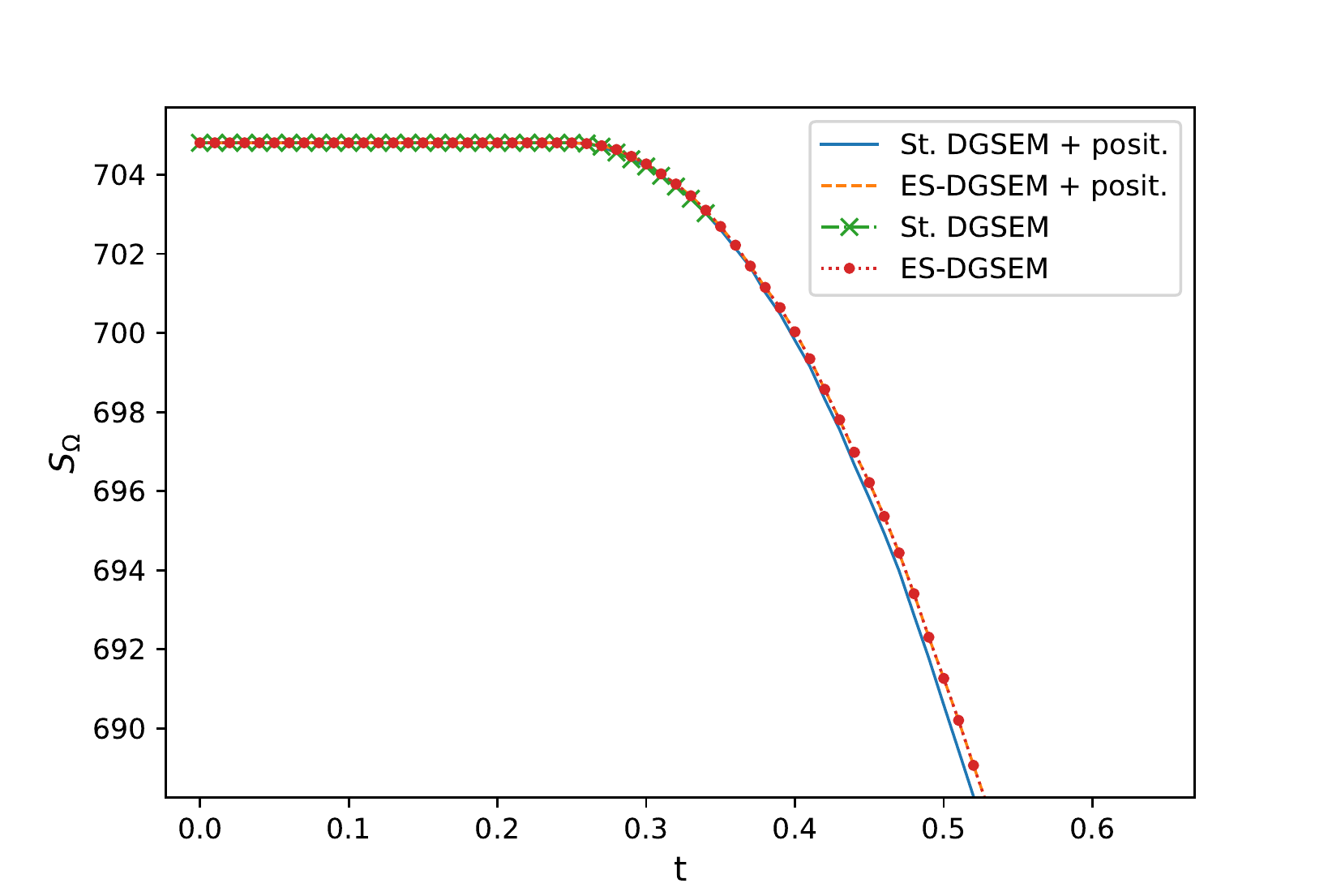}
}
\caption{Total entropy evolution over time (a) and detailed zoom at the beginning of the simulation (b) for the Sedov blast simulation.}
\label{fig:Sedov_Entropy}
\end{figure}

Figure \ref{fig:Sedov_Entropy} shows the evolution of the total mathematical entropy over time in the entire simulation domain.
It can be observed that, in spite of the use of the shock indicator, both the standard and the ES-DGSEM need positivity-preserving stabilization.
The standard DGSEM without FV stabilization crashes due to positivity problems when the shocks start forming at $t \approx 0.34$.
The ES-DGSEM simulation without stabilization is able to run longer, but finally crashes 
at $t \approx 12.5$.
However, both the standard and the ES-DGSEM are able to run until the end of the simulation with the proposed positivity limiter strategy.

\begin{figure}[htb]
\centering
\subfigure[Maximum blending coefficient.]{
	\includegraphics[trim=20 0 40 0,clip,width=0.42\linewidth]{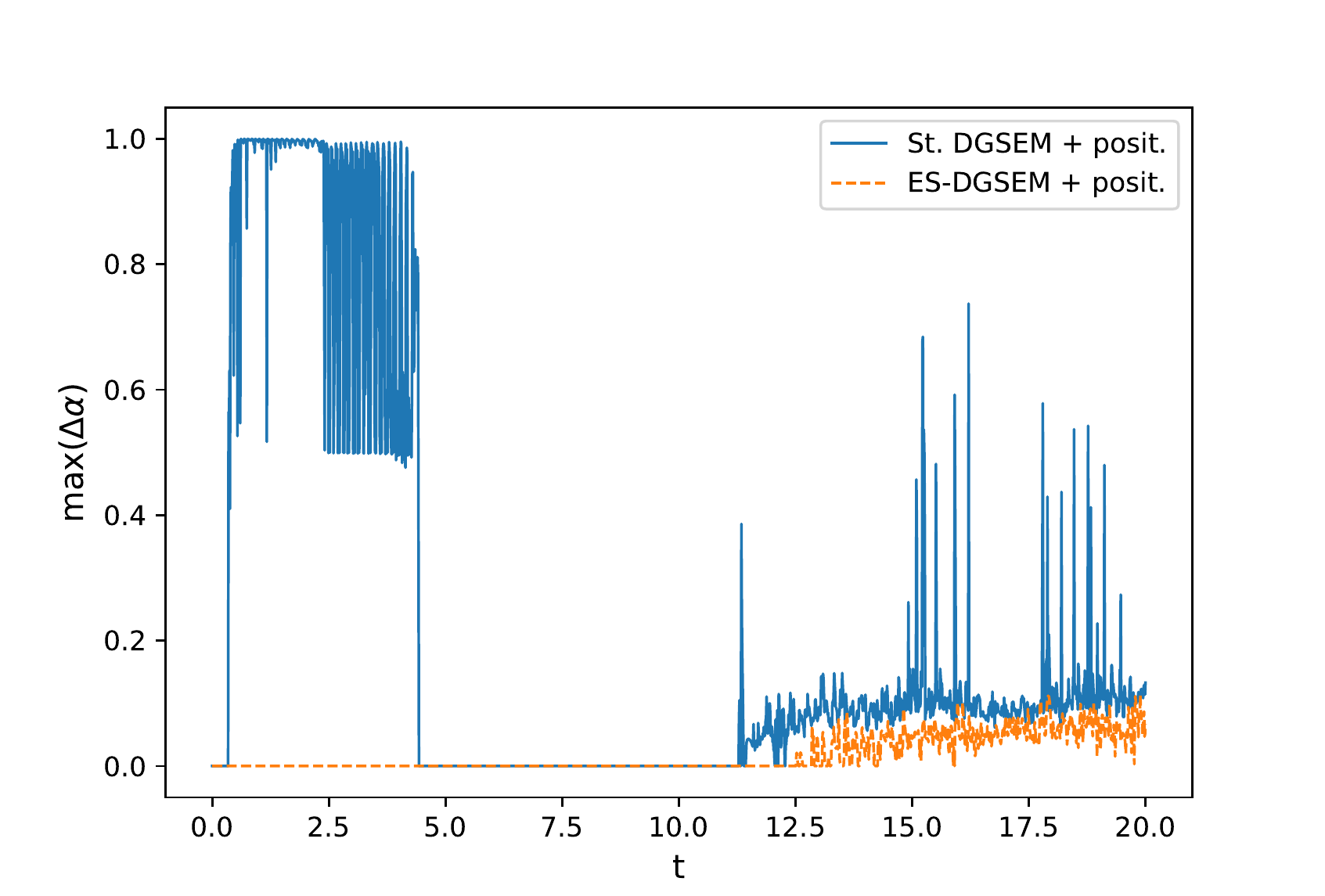}
}
\subfigure[Amount of low order method used in the domain.]{
	\includegraphics[trim=20 0 40 0,clip,width=0.42\linewidth]{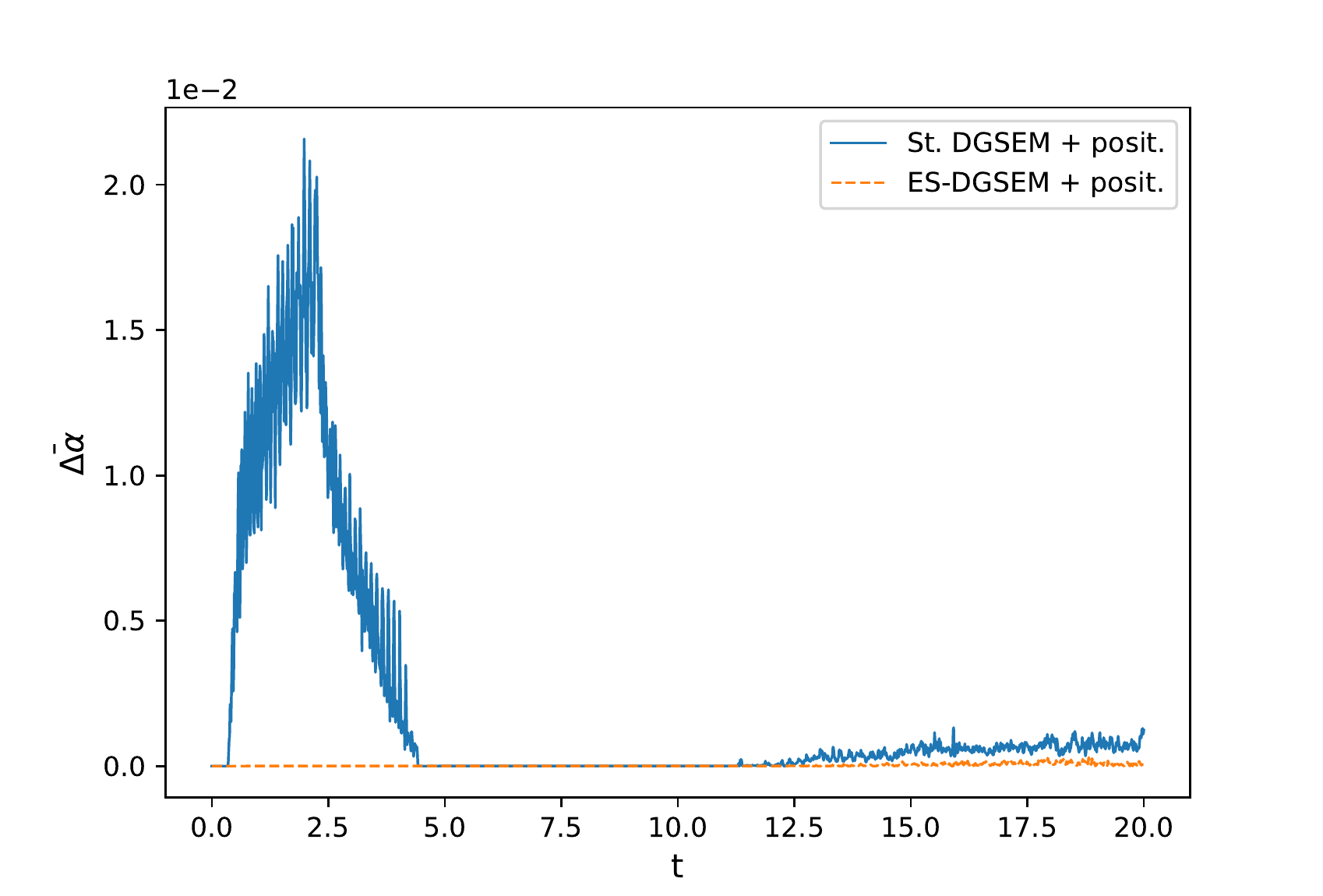}
}
\caption{Evolution of the \textit{correction} to the blending coefficient over time for the Sedov blast simulation.}
\label{fig:SedovBlendingCoeff}
\end{figure}

Figure \ref{fig:SedovBlendingCoeff} shows the evolution of the correction to the blending coefficient over time for the entire simulation using a sample time $\Delta \tau = 0.01$.
Note that the illustrated quantities are computed with \eqref{eq:maxAlpha} using $\Delta \alpha = \alpha_{\mathrm{new}} - \alpha$ instead of $\alpha$, such that we can appreciate the action of the positivity limiter.

As can be observed in Figure \ref{fig:SedovBlendingCoeff}(a), the ES-DGSEM only needs a small amount of extra first-order dissipation to stabilize the last part of the simulation ($t>12.5$), where the blending coefficient correction is always below $\Delta \alpha \le 0.11$.
On the other hand, the standard DGSEM needs very high corrections of the blending coefficient at the beginning of the simulation, where $\Delta \alpha$ oscillates between $0.5$ and $1$.
It is clear that the shock indicator, \eqref{eq:ShockIndicator}, was tuned for the ES-DGSEM.

Even though the standard DGSEM needs high amounts of stabilization in some elements, Figure \ref{fig:SedovBlendingCoeff}(b) shows that the amount of the domain that is discretized with a first-order FV correction is less than $2.15 \%$ at the beginning of the simulation, and less than $0.13\%$ at the end of the simulation when the base-line scheme is the standard DGSEM.
Moreover, when the baseline scheme is the ES-DGSEM, the amount of first-order FV correction is less than $0.03 \%$.

\begin{figure}[ht]
\centering
\includegraphics[trim=700 1420 700 30,clip,width=0.48\linewidth]{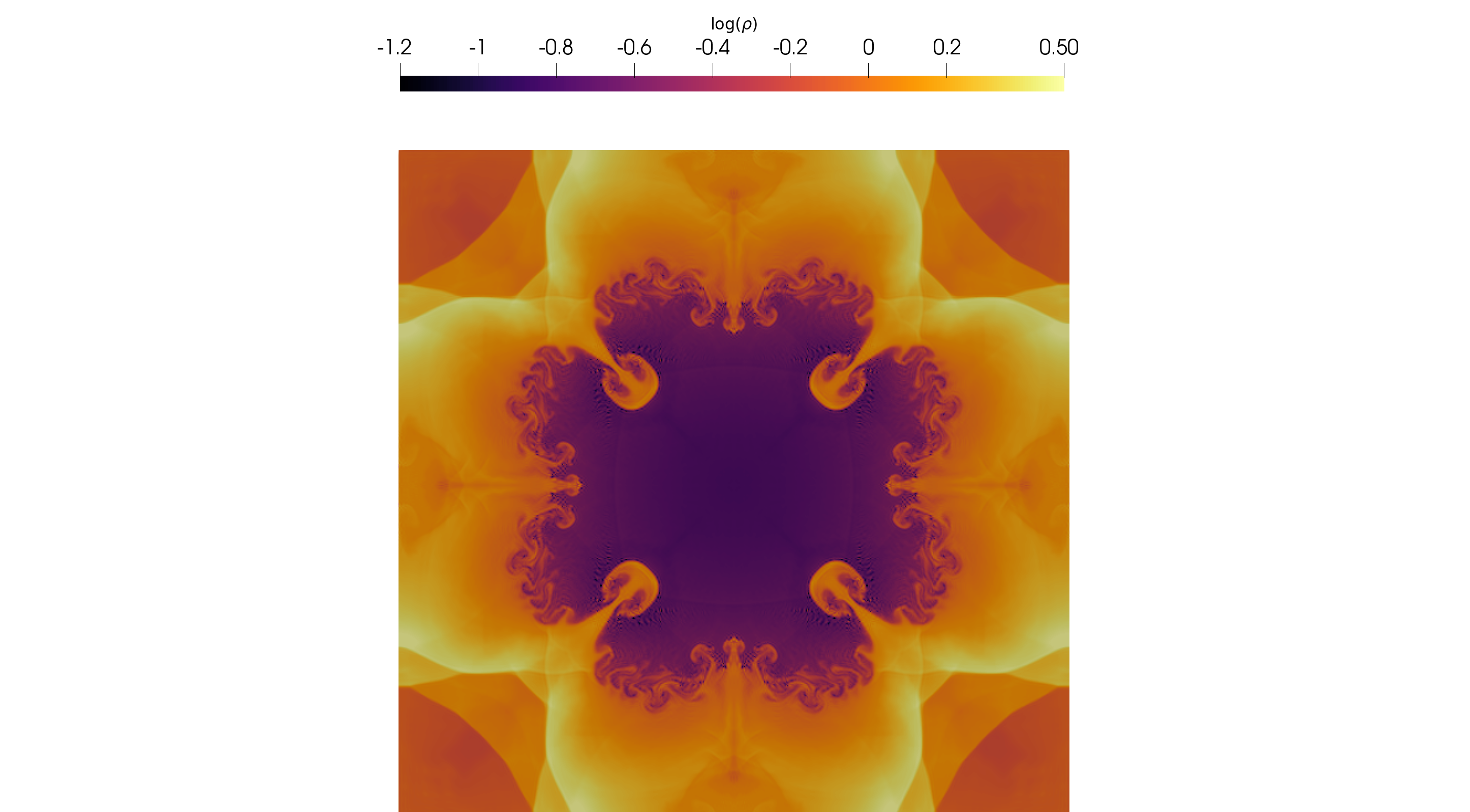}

\subfigure[St. DGSEM, $t=2$.]{
	\includegraphics[trim=700 100 700 100,clip,width=0.33\linewidth]{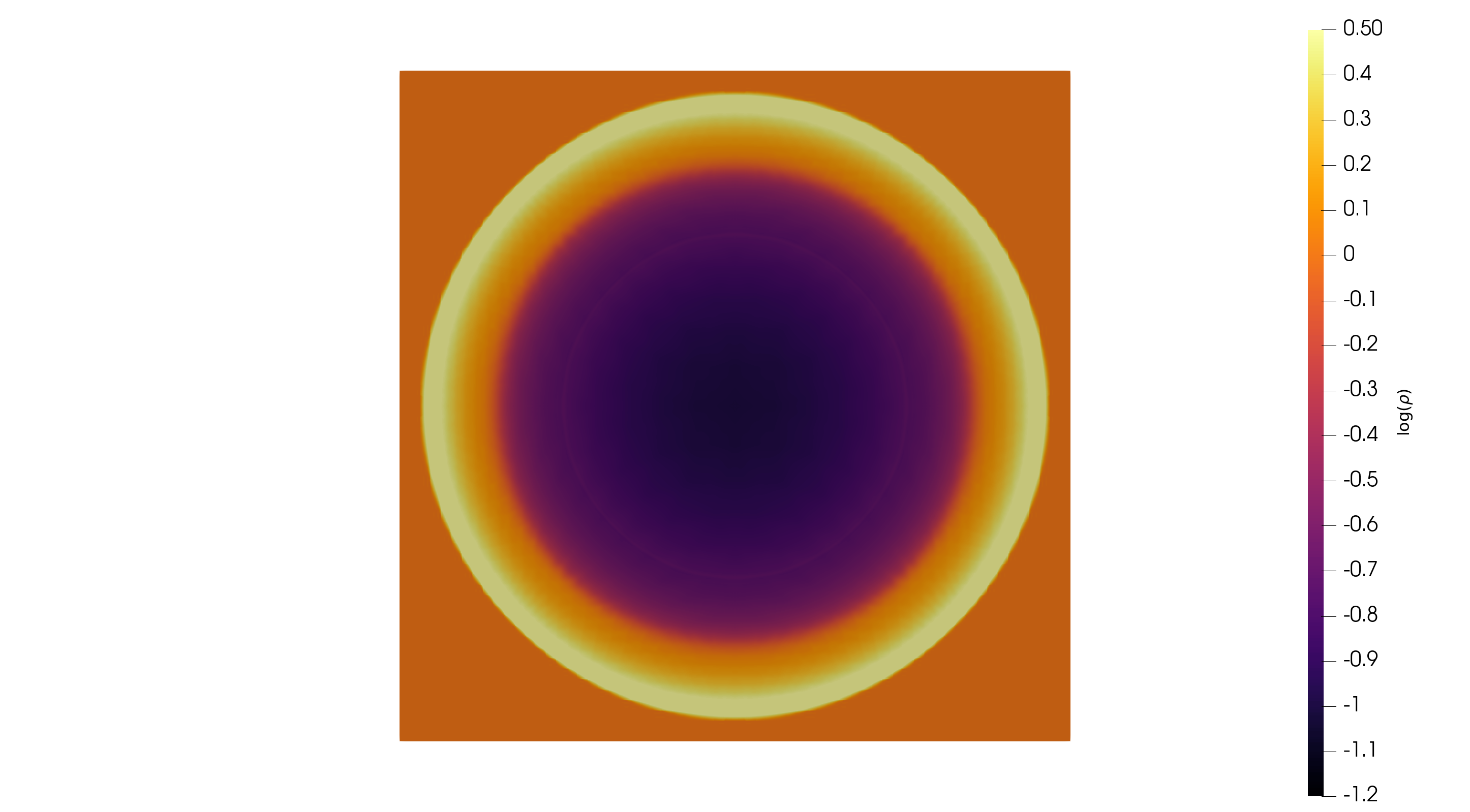}
}
\subfigure[St. DGSEM, $t=12.6$.]{
\includegraphics[trim=700 100 700 100,clip,width=0.33\linewidth]{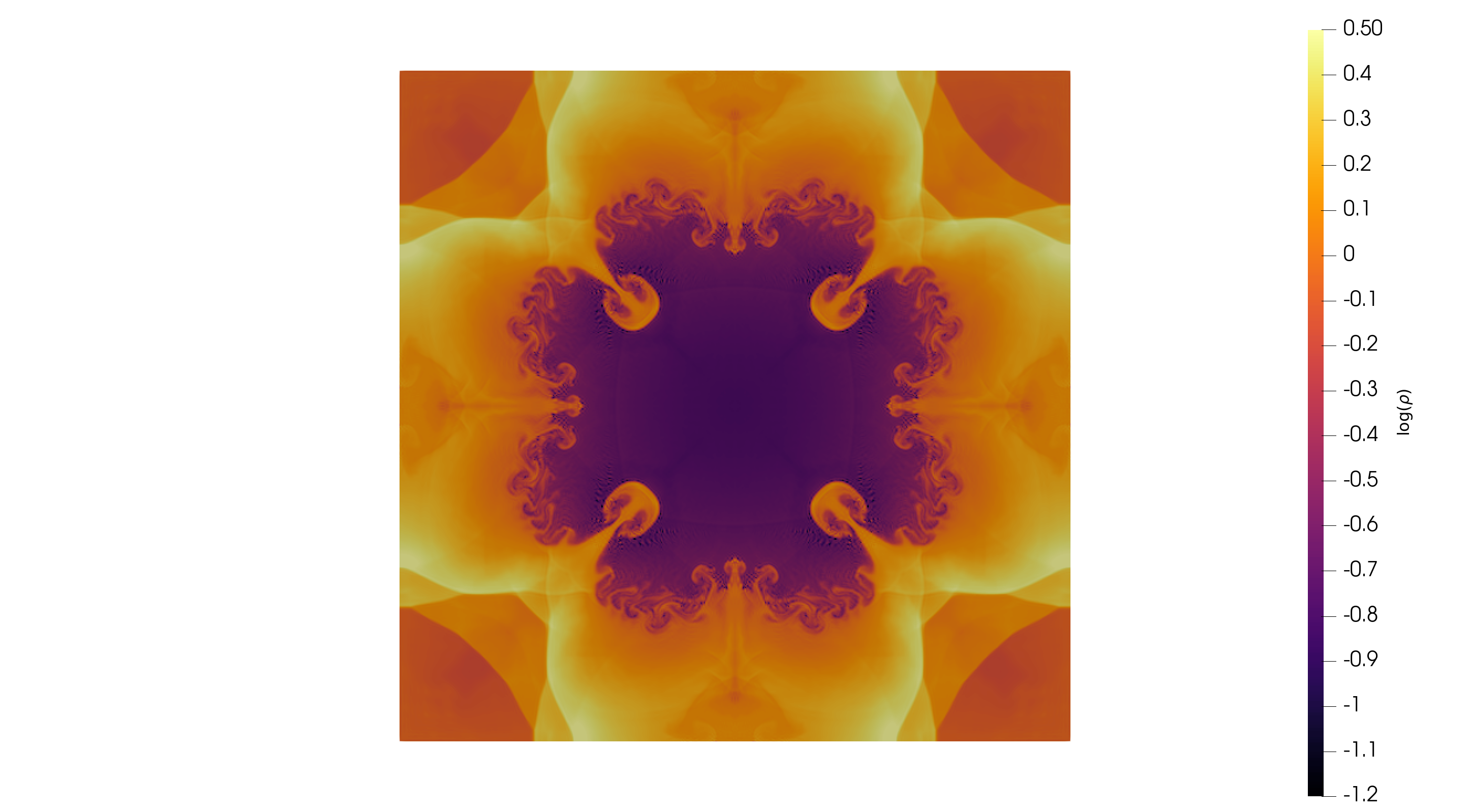}
}

\subfigure[ES-DGSEM, $t=2$.]{
	\includegraphics[trim=700 100 700 100,clip,width=0.33\linewidth]{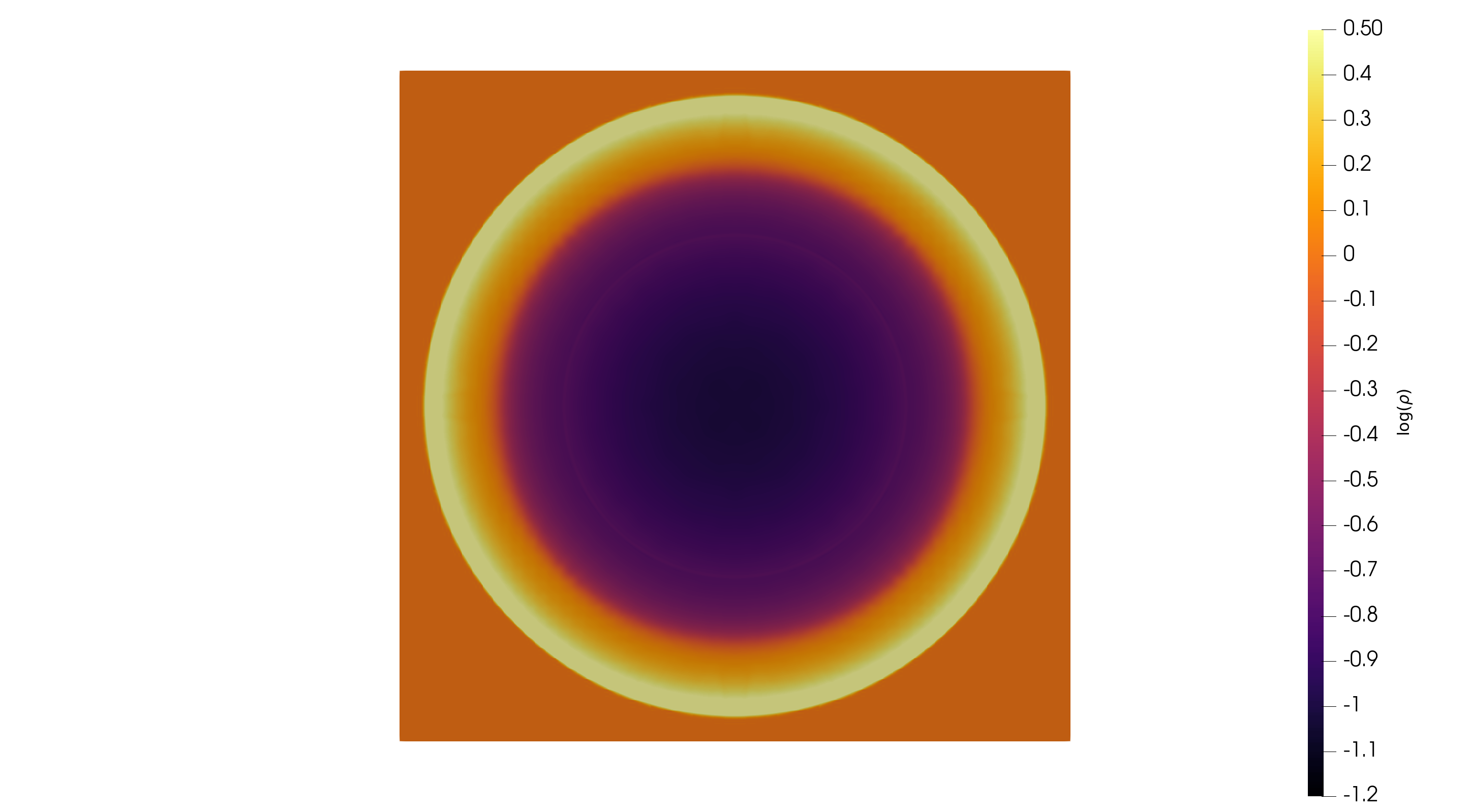}
}
\subfigure[ES-DGSEM, $t=12.6$.]{
	\includegraphics[trim=700 100 700 100,clip,width=0.33\linewidth]{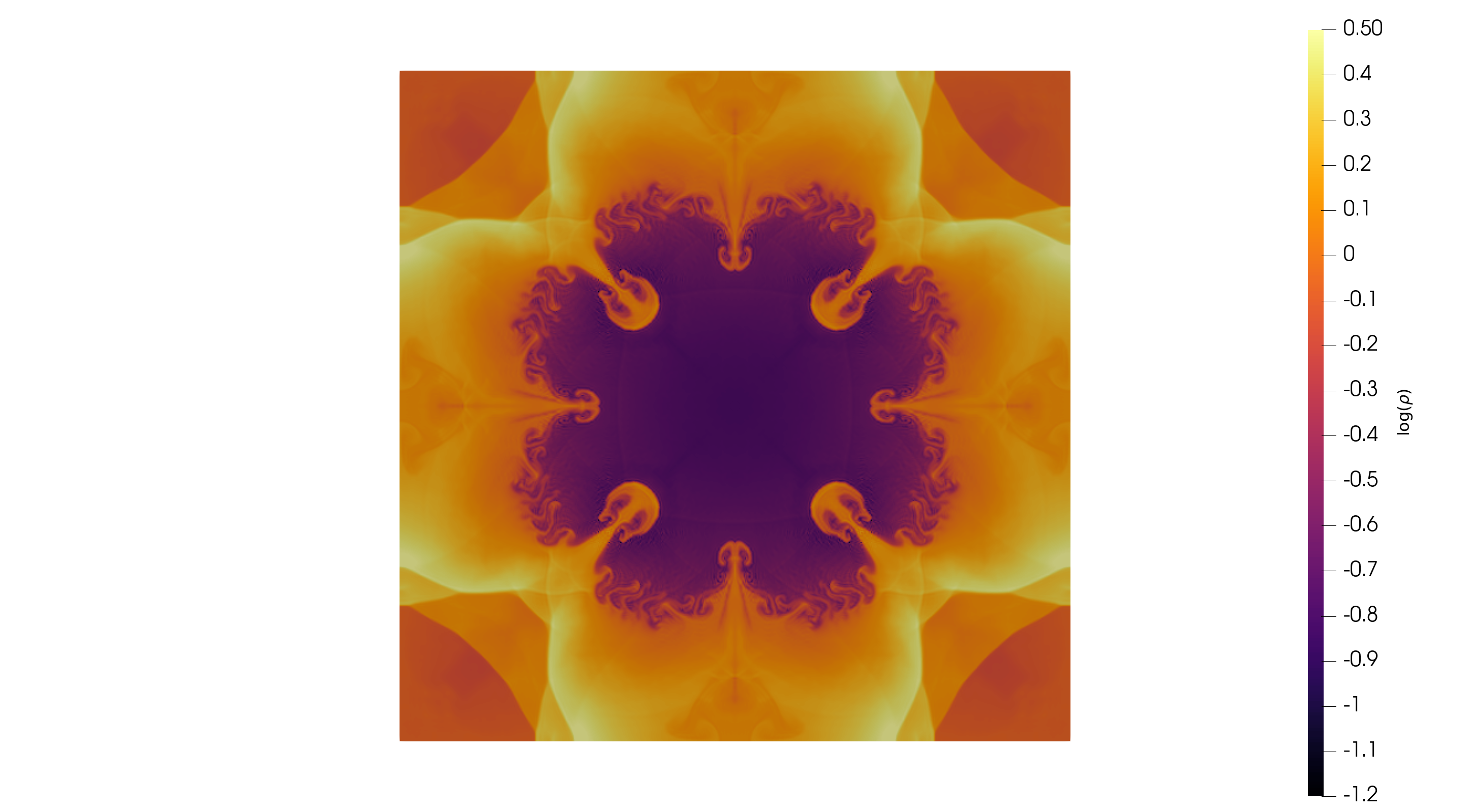}
}

\caption{Snapshots of the density logarithm for the Sedov blast simulations at $t=2$ (left) and $t=12.6$ (right) for the standard DGSEM (top) and the ES-DGSEM (bottom).}
\label{fig:Sedov_ContoursDens}
\end{figure}

In Figure \ref{fig:Sedov_ContoursDens}, we 
plot $\log(\rho)$ for different snapshots, as this quantity allows to appreciate the different scales that appear in the problem.
Note that the standard DGSEM solutions are more distorted than the ones obtained with the ES-DGSEM, as was also observed in the KHI simulation (Figure \ref{fig:KHI_Contours}).

Vortices develop in the domain due to the appearance of Kelvin-Helmholtz instabilities in the density layer, which are in turn triggered by the shock-shock interactions.
In fact, it is the shock-driven vortical dominated flow part that causes the ES-DGSEM without positivity limiter to crash (see Figures \ref{fig:Sedov_Entropy}(a), \ref{fig:Sedov_ContoursDens}(b) and \ref{fig:Sedov_ContoursDens}(d)).
On the other hand, Figures \ref{fig:Sedov_ContoursDens}(a) and \ref{fig:Sedov_ContoursDens}(c) confirm that the standard DGSEM without positivity limiter crashes when shocks are forming, a long time before the shock front bounces on the boundaries.

\begin{figure}[ht]
\centering
\includegraphics[trim=700 1420 700 30,clip,width=0.48\linewidth]{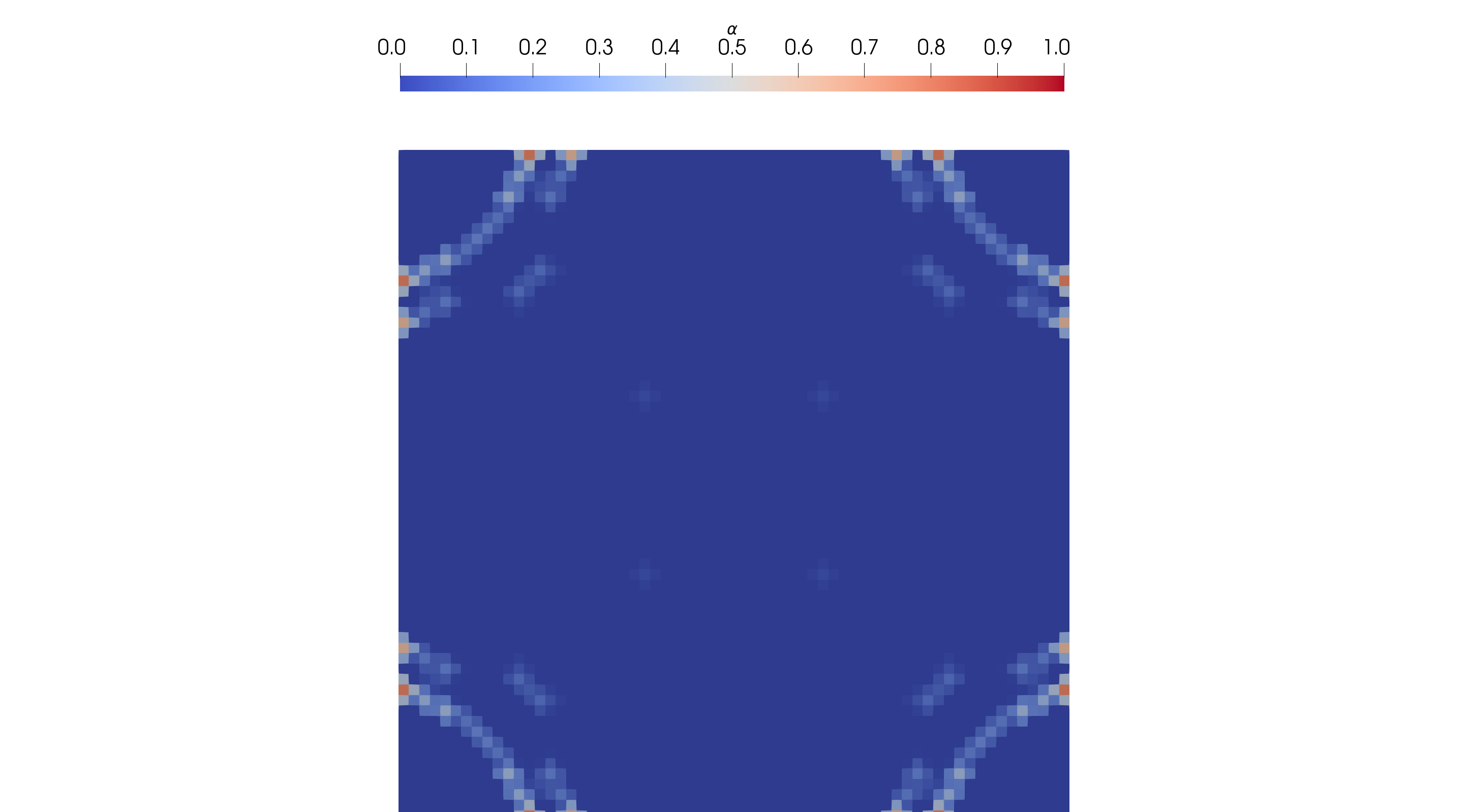}

\subfigure[St. DGSEM, $t=2$.]{
	\includegraphics[trim=700 100 700 100,clip,width=0.33\linewidth]{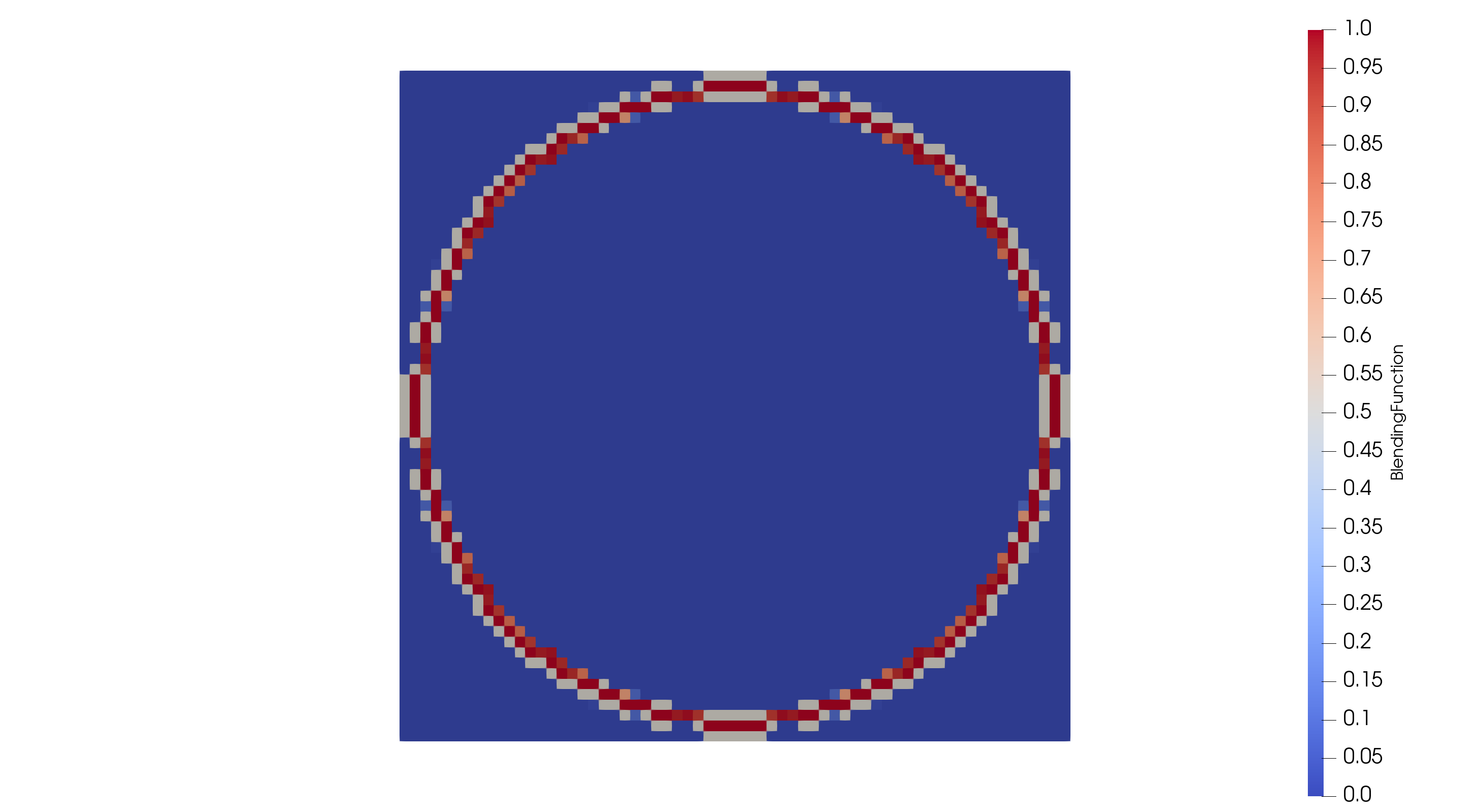}
}
\subfigure[St. DGSEM, $t=12.6$.]{
	\includegraphics[trim=700 100 700 100,clip,width=0.33\linewidth]{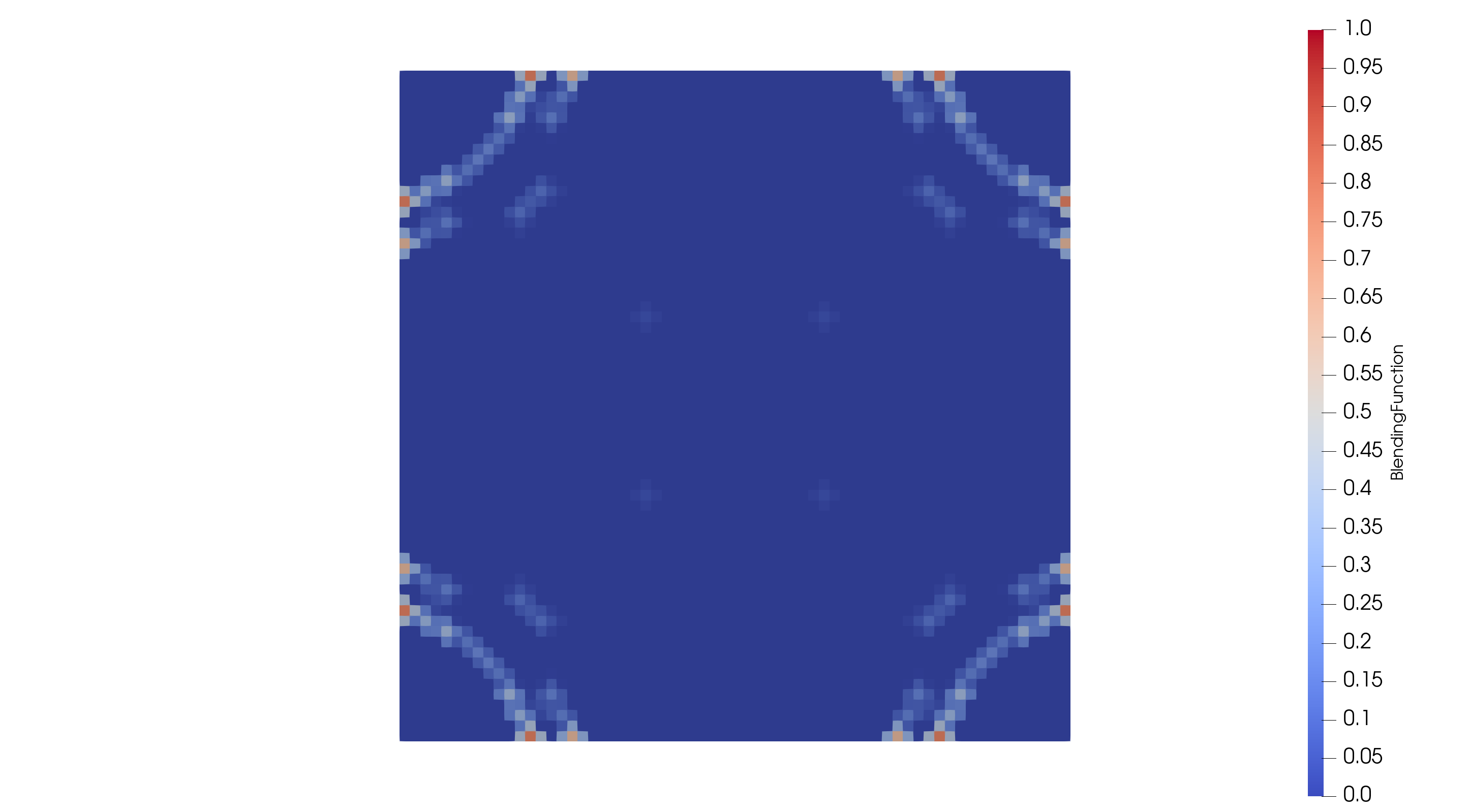}
}
\subfigure[ES-DGSEM, $t=2$.]{
	\includegraphics[trim=700 100 700 100,clip,width=0.33\linewidth]{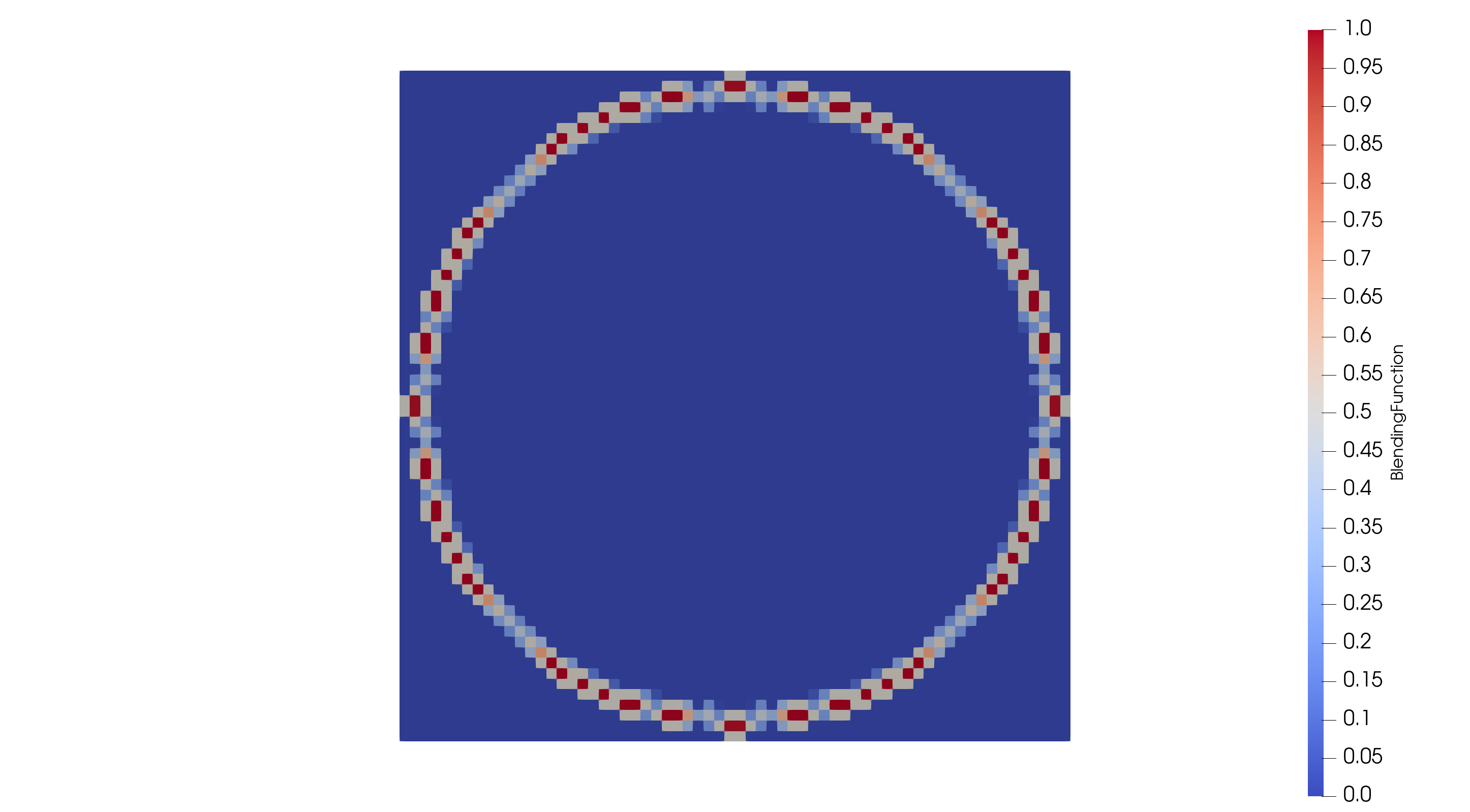}
}
\subfigure[ES-DGSEM, $t=12.6$.]{
	\includegraphics[trim=700 100 700 100,clip,width=0.33\linewidth]{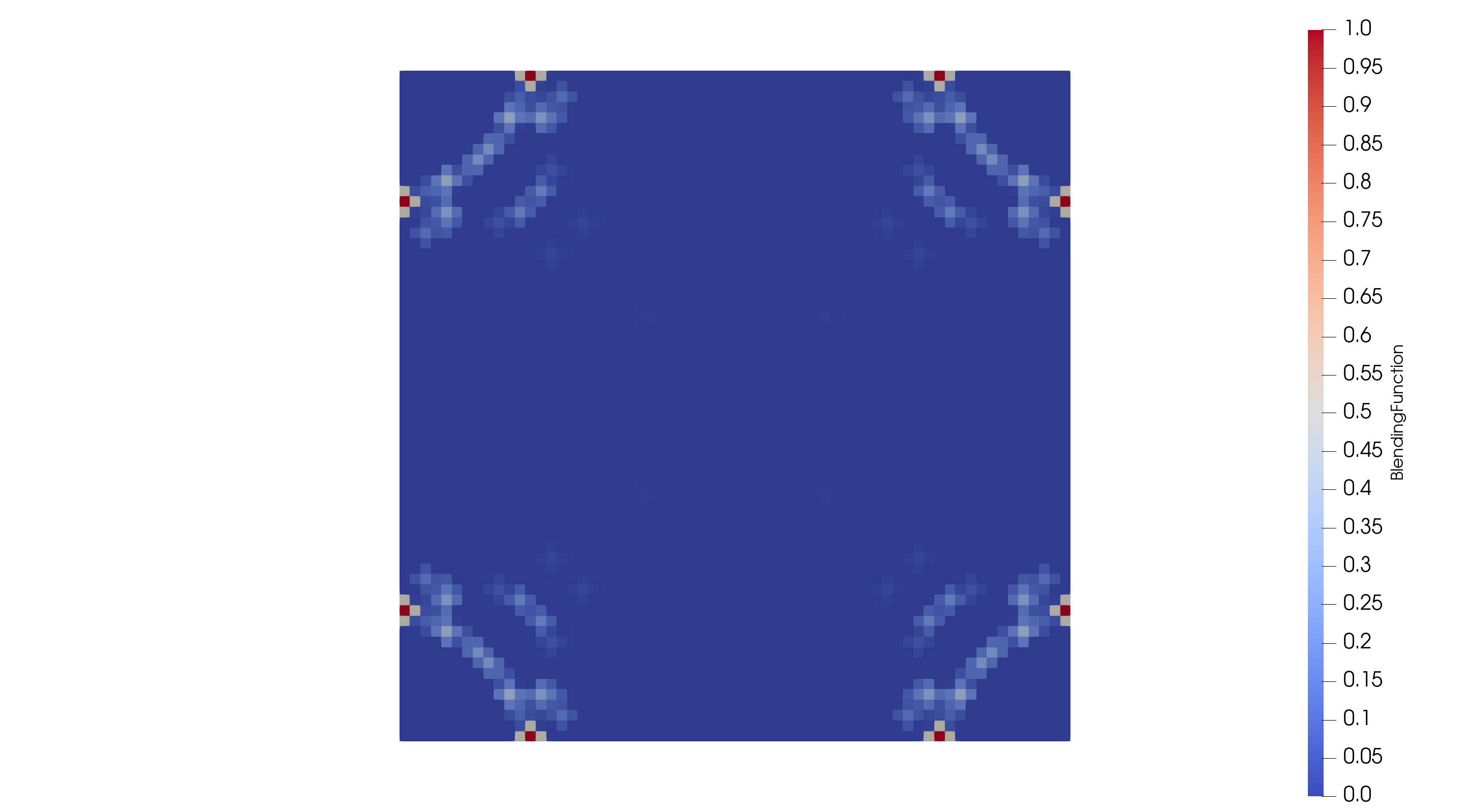}
}

\caption{Snapshots of the blending coefficient for the Sedov blast simulations at $t=2$ (left) and $t=12.6$ (right) for the standard DGSEM (top) and the ES-DGSEM (bottom).}
\label{fig:Sedov_ContoursAlpha}
\end{figure}

Finally, Figure \ref{fig:Sedov_ContoursAlpha} shows the distribution of $\alpha$ for snapshots of the simulations at $t=2$ and $t=12.6$.
A comparison of the blending coefficient contours with Figure \ref{fig:Sedov_ContoursDens} reveals that the indicator, \eqref{eq:ShockIndicator}, correctly detects the shocks and is not triggered on turbulent areas of the domain.
The action of the positivity limiter is evident in Figure \ref{fig:Sedov_ContoursAlpha}(a): the blending coefficient correction affects elements in the post-shock region of the standard DGSEM, where the spatial propagation of the blending coefficient had already acted.
This behavior reveals that an appropriate shock indicator for the standard DGSEM should detect the post-shock regions and explains the jumps of $\Delta \alpha$ between $0.5$ and $1$ that we observed in Figure \ref{fig:SedovBlendingCoeff}(a).

\section{Conclusions}
We presented a positivity-preserving limiter for DGSEM discretizations of the compressible Euler equations that is based on a subcell Finite Volume method.
We show that our strategy is able to ensure robust simulations with positive density and pressure when using the standard and split-form DGSEM if a positivity preserving Riemann solver is used.
We tested our scheme with the Rusanov and HLLE Riemann solvers in extremely under-resolved vortex dominated simulations and in problems with shocks.

\textbf{Acknowledgments} This work has received funding from the European Research Council through the ERC Starting Grant EXTREME, with the ERC grant agreement no. 714487.

\bibliographystyle{plainnatnourl}
\bibliography{Biblio.bib}

\end{document}